\documentclass[aos, authoryear]{imsart}
\RequirePackage[OT1]{fontenc}
\RequirePackage{amsthm,amsmath,natbib, comment}
\RequirePackage{hypernat}
\startlocaldefs
\numberwithin{equation}{section}
\theoremstyle{plain}

\endlocaldefs
\usepackage{psfrag, epsfig, graphicx}
\usepackage{amsfonts,amsmath,latexsym,amssymb}
\usepackage[active]{srcltx}
\usepackage{geometry}
\setattribute{journal}{name}{ }

\startlocaldefs
\numberwithin{equation}{section}
\theoremstyle{plain}
\newcommand{\epr}{\hfill\hbox{\hskip 4pt
                \vrule width 5pt height 6pt depth 1.5pt}\vspace{0.5cm}\par}
\renewcommand{\(}{$\,}
\renewcommand{\)}{\,$}

\newcommand{\cc}[1]{\mathcal{#1}}

\renewcommand{\hat}[1]{\widehat{#1}}
\renewcommand{\tilde}[1]{\widetilde{#1}}
\newcommand{\nn}{\nonumber \\}
\numberwithin{equation}{section}
\numberwithin{figure}{section}
\newtheorem{assumption}{Assumption}
\newtheorem{theorem}{Theorem}
\newtheorem{proposition}{Proposition}
\newtheorem{definition}{Definition}
\newtheorem{lemma}{Lemma}

\def\T{\top}

\def\dd{\mathrm{d}}
\def\RR{\mathbb{R}}

\def\RRd{\mathbb{R}^{d}}
\def\SSd{\mathbb{S}^{d-1}}
\def\SS{\mathbb{S}^{1}}

\def\supp{\operatorname{supp}}

\def\kappa{\varkappa}
\def\ta{\theta}
\def\eps{\varepsilon}

\def\Var{\operatorname{Var}}
\def\EE{\mathbb{E}}
\def\PP{\mathbb{P}}
\def\TH{\operatorname{TH}}

\newcommand{\cQ}{{\cal Q}}

\newcommand{\expect}[1]{\EE^{\eps}_{#1}}
\newcommand{\norm}[2]{\cc{N} \left( {#1},{#2} \right)}
\def\bob{\beta_{\max}}

\def\numhyp{N_{\eps}}

\newcommand{\ort}[1]{{#1}_{\bot}}
\def\C(S){1}                                        %% constant bounding the bias
\def\K{\cc K}                                       %% kernel function
\endlocaldefs

\begin{document}

\begin{frontmatter}
\title{Adaptive estimation in the single-index model\\via oracle approach \thanksref{T1}}
\runtitle{Adaptation in the single-index model}
\thankstext{T1}{Funding of the ANR-07-BLAN-0234 is acknowledged. The second author is also supported by the DFG FOR 916.}

\begin{aug}
\author{\fnms{Oleg} \snm{Lepski}%\thanksref{t1,t2,m1}
\ead[label=e1]{ Oleg.Lepski@cmi.univ-mrs.fr}}%,
\and \author{\fnms{Nora} \snm{Serdyukova}%\thanksref{t3,m1,m2}
\ead[label=e3]{Nora.Serdyukova@gmail.com}}
%\and
%\author{\fnms{Third} \snm{Author}\thanksref{t1,m2}
%\ead[label=e3]{third@somewhere.com}
%\ead[label=u1,url]{http://www.foo.com}}

%\thankstext{t1}{Some comment}
%\thankstext{t2}{First supporter of the project}
%\thankstext{t3}{Second supporter of the project}
\runauthor{O. Lepski and N. Serdyukova}

\affiliation{Aix-Marseille Universit\'e %\thanksmark{m1}
   and Georg-August-Universit\"at G\"ottingen%\thanksmark{m2}
}

\address{Laboratoire d'Analyse, \\Topologie, Probabilit\'es UMR 7353\\
Aix-Marseille Universit\'e \\39, rue F. Joliot Curie \\13453 Marseille FRANCE\\
%Usually a few lines long\\
\printead{e1}\\
\phantom{E-mail:\ }%\printead*{e2}
}

\address{Institute for Mathematical Stochastics \\Georg-August-Universit\"at G\"ottingen \\Goldschmidtstra{\ss}e 7 \\37077 G\"ottingen GERMANY\\
%Usually a few lines long\\
%Usually a few lines long\\
\printead{e3}\\
%\printead{u1}
}
\end{aug}

\begin{abstract}
In the framework of nonparametric multivariate function estimation we are interested in structural adaptation. We assume that the function to be estimated has the ``single-index'' structure where neither the link function nor the index vector is known. We suggest a novel procedure that adapts simultaneously to the unknown index and smoothness of the link function. For the proposed procedure, we prove a ``local'' oracle inequality (described by the pointwise seminorm), which is then used to obtain the upper bound on the maximal risk of the adaptive estimator under
assumption that the link function belongs to a scale of H\"{o}lder classes.
The lower bound on the minimax risk shows that in the case of estimating at a given point the constructed estimator is optimally rate adaptive over the considered range of classes. For the same procedure we also establish a ``global'' oracle inequality  (under the $ L_r $ norm, $r< \infty $) and examine its performance over the Nikol'skii classes. This study shows that the proposed method can be applied to estimating functions of inhomogeneous smoothness, that is whose smoothness may vary from point to point.

%The problem of adaptive estimation of a multivariate function under the single-index constrains when
%both the link function and index vector are unknown is investigated. We propose a novel estimation procedure that adapts simultaneously
 %to the unknown index vector and the smoothness of link function by selecting from a family of specific kernel
 %estimators. We establish a pointwise oracle inequality which, in its turn, is used to judge the quality
 %of estimating the entire function as well (global oracle inequality). Both results are applied to  the problems of pointwise and global
 %adaptive estimation  over a collection of H\"older and Nikol'skii functional classes.

\end{abstract}



\end{frontmatter}

%\begin{thm}
%All conjectures are interesting, but some conjectures are more
%interesting than others.
%\end{thm}
%%%%%%%%%%%%%%%%%%%%%%%%%%%%%%%%%%%%%%%%%%%%%%%%%%%%%%%%%%%%%%%%%%%%%%%%%%%%%%%%%%%%%%%%%%%%%%%%%%%%
%%%%%%%%%%%%%%%%%%%%%%%%%%%%%%%%%%%%%%%%%%%%%%%%%%%%%%%%%%%%%%%%%%%%%%%%%%%%%%%%%%%%%%%%%%%%%%%%%%%%
\section{Introduction}\label{Section_Intro}
%%%%%%%%%%%%%%%%%%%%%%%%%%%%%%%%%%%%%%%%%%%%%%%%%%%%%%%%%%%%%%%%%%%%%%%%%%%%%%%%%%%%%%%%%%%%%%%%%%%%
%%%%%%%%%%%%%%%%%%%%%%%%%%%%%%%%%%%%%%%%%%%%%%%%%%%%%%%%%%%%%%%%%%%%%%%%%%%%%%%%%%%%%%%%%%%%%%%%%%%%
This research aims at estimating multivariate functions with the use of the oracle approach. The first step of the method consists in justification of pointwise and global oracle inequalities for the estimation procedure; the second step is the deriving from them adaptive results for estimation of the point functional and the entire function correspondingly. The obtained results show full adaptivity of the proposed estimator as well as its minimax rate optimality.

%This paper deals with  estimation of  multivariate functions. We establish local as well as global oracle inequalities and show how to use them for deriving minimax adaptive results.
\paragraph{ Model and set-up}
Let \( \cc D \supset [-1/2, 1/2]^d \) be a bounded interval in \( \RRd \). We observe  a path \( \{ Y_{\eps}(t), t \in \cc D \} \), satisfying the stochastic differential equation
\begin{equation}\label{WGN_model}
    Y_{\eps}(\dd t) = F(t) \dd t + \eps W(\dd t) \; , \;\; t = (t_1, \ldots, t_d) \in \cc D,
\end{equation}
where \( W \) is a Brownian sheet and \( \eps \in (0,1) \) is the  deviation parameter.

\par In the single-index modeling the signal \( F \) has a particular structure:
\begin{equation}\label{single-index}
    F(x) = f(x^{\T} \ta^{\circ}),
\end{equation}
where \( f : \RR \to \RR \) is called link function and \( \ta^{\circ} \in \SSd \) is the index vector.

\par %We consider the observation domain \( \cc D \) larger than \( [-1/2, 1/2]^d \) in order to skip the study of boundary effects.
%Moreover, to simplify the presentation of our approach and results
We consider the case of completely unknown parameters \( f \) and \( \ta^{\circ} \) and the only technical assumption is that $f\in\mathbb{F}_M$ where
$\mathbb{F}_M=\left\{g: \RR \to \RR\; |\; \sup_{u\in\RR}|g(u)| \le M\right\}$ for some $M>0$. However, the knowledge of $M$ as well as any information on the smoothness of the link function are not required for the proposed below estimation procedure. The consideration is restricted to the case $d=2$ except the second assertion of Theorem \ref{th:pointwise-adaptation} concerning a lower bound for function estimation at a given point.
Also, without loss of generality we will assume that $\cc D=[-1,1]^{2}$ and $\eps\le e^{-1}$.

\par Let \( \tilde F(\cdot) \) be an estimator, i.e. a measurable function of the observation
\( \{ Y_{\eps}(t), t \in \cc D \} \) and \( \EE_F^{\eps} \) denote the mathematical expectation with respect to \( \PP_F^{\eps} \), the family of probability distributions generated by the observation process \( \{ Y_{\eps}(t), t \in \cc D \} \) on the Banach space of continuous functions on \( \cc D \), when \( F \) is the mean function. The estimation quality is measured by  the \( L_r \) risk, \( r \in [1, \infty) \),
\begin{equation}\label{global risk}
    \cc R_r^{(\eps)} (\tilde F, F)
        =
    \expect{F} \| \tilde F  - F \|_r,
\end{equation}
where \( \| \cdot \|_r \) is the \( L_r \) norm on  \( [-1/2, 1/2]^2 \)
or by the ``pointwise''  risk
\begin{equation}\label{local risk}
    \cc R_{r,x}^{(\eps)} (\tilde F, F)
        =
    (\expect{F} |\tilde F (x) - F(x)|^r )^{1/r}.
\end{equation}

%Also for the clarity of presentation we will assume
%that $\ta^{\circ}\in\SS$, where $\SS$ is the unite circle in \( \RR^2 \) Section~\ref{sec:extensions } it is shown that our results can be extended to the case $\ta^{\circ}\in\mathbb{R}^2$.

\par The aim is to estimate the entire function \( F \) on
\( [-1/2, 1/2]^2 \) or its value \( F(x) \) from the observation \( \{ Y_{\eps}(t), t \in \cc D \} \) satisfying SDE
\eqref{WGN_model} without any prior knowledge of the nuisance parameters: the function \( f \) and the unit vector \( \ta^{\circ} \). More precisely, we will construct an adaptive (not depending of \( f \) and \( \ta^{\circ} \)) estimator \( \hat F(x) \) at any point \( x \in [-1/2,1/2]^{2} \). In what follows \( \hat F \) notation stands for an adaptive estimator and \( \tilde F \) denotes an arbitrary estimator. Our estimation procedure is a random selector from a special family of kernel estimators parametrized by a window size (bandwidth) \( h>0 \) and a direction of the projection \( \ta \in \SS \), see Section \ref{sect:SelectionRule} below. For that procedure we then establish a pointwise oracle inequality (Theorem \ref{th:local-oracle-inequality}) of the following type:
\begin{equation}\label{eq:intro_pointwise_oracle_ineq_rough}
    \cc R_{r,x}^{(\eps)} (\hat F, F)
    \le
    C_1 \, \eps \sqrt{\ln(1/\eps)  /h^{*} (x^{\T} \ta^{\circ})}
        + C_2  \, \eps \sqrt{\ln(1/\eps)  },
\end{equation}
where \( h^{*} \) is an optimal in a certain sense (oracle) bandwidths, see Definition~\ref{def_oracle bandwidth}. As
\( r < \infty \) Jensen's inequality and Fubini's theorem trivially imply
\begin{equation*}
    \left[ \cc R_r^{(\eps)} (\hat F, F)\right]^r
    \le
    \expect{F} \big\| \hat F (\cdot) - F(\cdot) \big\|^{r}_r
    =
    \big\|\cc R_{r,\cdot}^{(\eps)} (\hat F, F)\big\|_r^r.
\end{equation*}
Hence, we immediately obtain the ``global'' oracle inequality
\begin{equation}\label{eq:intro_global_oracle_ineq_rough}
    \cc R_{r}^{(\eps)} (\hat F, F)
    \le
    C_1 \, \eps \| \sqrt{\ln(1/\eps)  /h^{*}} \|_r
    + C_2  \, \eps \sqrt{\ln(1/\eps)  }.
\end{equation}
Both inequalities \eqref{eq:intro_pointwise_oracle_ineq_rough} and \eqref{eq:intro_global_oracle_ineq_rough} aside of being quite informative itself -- we will see in Section~\ref{OracleEstmator} from Proposition~\ref{prop:risk-of-oracle-estimator} that they claim that our adaptive estimator mimics its ideal (oracle) counterpart, i.e. their risk bounds differ only by a numerical constant, -- they are further used to judge the minimax rate of convergence under the pointwise and \( L_r \) losses correspondingly (Theorems \ref{th:pointwise-adaptation} and \ref{th:global-adaptation}). We will see that these rates are in accordance with Stone's dimensionality reduction principle, see pp. 692-693 in \cite{Stone85}. Indeed, as the statistical model is effectively one-dimensional due to the structural assumption \eqref{single-index} so the rate of convergence is.

\par The obtained results demonstrate full adaptivity of the proposed estimator to the unknown direction of the projection \( \ta^{\circ} \) and the smoothness of \( f \). Moreover, the lower bound given in the second assertion of Theorem~\ref{th:pointwise-adaptation} shows that in the case of pointwise estimation over the range of classes of \( d \)-variate functions having the single-index structure, see definition~\eqref{Def:SI-class_local}, our estimator is even optimally rate adaptive, that is it achieves the minimax rate of convergence. This fact is in striking contrast to the common knowledge that a payment for pointwise adaptation in terms of convergence rate is unavoidable. Indeed, if the index \( \ta^{\circ} \) would be known, than the problem boils down to pointwise adaptation over H\"{o}lder classes in the univariate GWN model. As demonstrated in \cite{Lep1990}, an optimally adaptive estimator does not exist in this case.

\par Although the literature on the single-index model is rather numerous, we mention only books \cite{Haerdle2004}, \cite{Horowitz1998}, \cite{Gyorfi2002}  and \cite{Koros}, quite a few works address the problem of function estimating when both the link function and index are unknown. To the best of our knowledge the only exceptions are \cite{Golubev1992}, \cite{Lecue} and \cite{GoldLep2008}. An adaptive projection estimator is constructed in \cite{Golubev1992}, in \cite{Lecue} the aggregation method is used. Both the papers employ \( L_2 \) losses. \cite{GoldLep2008} seems to be the first work on pointwise adaptive estimation in the considered set-up, the upper bound for estimation at a point obtained therein is similar to our, but the estimation procedure is different.

%However, this approach requires the knowledge of the minimal smoothness of the link function and the bound on its \( sup \)-norm not only to obtain theoretical results but also as an input of the procedure.

\paragraph{Organization of the paper}  In Section \ref{subsec:oracle-approach} we motivate and explain the proposed
selection rule. Then in Section~\ref{subsec:OI} we establish for it local and global oracle inequalities of type \eqref{eq:intro_pointwise_oracle_ineq_rough} and \eqref{eq:intro_global_oracle_ineq_rough}. In Section \ref{subsec:adaptive-estimation} we apply these results to minimax adaptive estimation. Particularly, Section~\ref{subsec:PointwiseAdaptation} is devoted to the upper bound and already discussed above lower bound for estimation over a range of H\"{o}lder classes. Section \ref{sect:Global_Adapt} addresses the ``global'' adaptation under the \( L_r \) losses and the estimator performance over the collection of classes of single-index functions with the link function in a Nikol'skii class, see Definition~\ref{def:nikolskii-class} and~\eqref{Def:SI-class_global}. That consideration allows to incorporate in analysis functions of inhomogeneous smoothness, that is those which can be very smooth on some parts of observation domain and irregular on the others. The proofs of the main results are given in Section~\ref{sec:proofs} and the proofs of technical lemmas are postponed until Appendix.

\begin{comment}
\section{Main results}
\label{sec:main-results}
Then, we apply these results to adaptive estimation over a collection of H\"older classes (pointwise estimation) and over a collection of
Nikol'skii classes (estimating the entire function with the accuracy of an estimator measured under the~$L_r$~risk).
\end{comment}

\section{Oracle approach}
\label{subsec:oracle-approach}
Below we define an ``ideal'' (oracle) estimator and describe our estimation procedure. Then we present local and global oracle inequalities demonstrating a nearly oracle performance of the proposed estimator.
%In this section we present our procedure and establish for it local and global oracle inequalities.

\par Denote by $\K:\RR \to \RR$ any function (kernel) that integrates to one, and define for any $z\in \RR$, $h\in(0,1]$ and any $f\in\mathbb{F}_M$
$$
\Delta_{\K,f}(h,z)
=
\sup_{\delta\le h}\left|  \frac{1}{\delta}  \int\K\Big( \frac{u-z}{\delta}  \Big)\big[f(u)-f(z)\big]\dd u\right|,
$$
%We note that $1/\delta\int\K\big([u-z]/\delta\big)f(u)\dd u$ (kernel smoother) can be understood as an approximation of the function $f$ at the point $z$.
%Thus, $\Delta_{\K,f}(h,z)$ is a
a monotonous approximation error of the kernel smoother
\( 1/\delta\int\K\big[ (u-z)  /\delta\big] f(u)\dd u \). In particular, if the function $f$ is uniformly continuous then\( \Delta_{\K,f}(h,z)\to 0 \) as \( h\to 0 \).

In what follows we assume that the kernel $\K$ obeys
\begin{assumption}
\label{ass:assumption-on-kernel}
%\begin{enumerate}\item
(1)\; $\text{supp}(\K)\subseteq[-1/2,1/2]$, $\int \K=1$, $\K$ is symmetric;

\medskip

\hskip2.4cm (2)\; there exists $Q>0$ such that

\smallskip

\hskip3.2cm$
\big|\K(u)-\K(v)\big|\leq Q|u-v|,\quad\forall u,v\in \RR.
$

%\end{enumerate}

\end{assumption}

\subsection{Oracle estimator}\label{OracleEstmator} For any $y\in \RR$ denote by
$$
\overline{\Delta}_{\K,f}(h,y)=\sup_{a>0} \frac{1}{2a}  \int_{y-a}^{y+a}\Delta_{\K,f}(h,z)\dd z,
$$
the Hardy-Littlewood maximal function of $\Delta_{\K,f}(h,\cdot)$, see for instance \cite{WheedenZygmund1977}.
Put also $\Delta^{*}_{\K,f}(h,\cdot)=\max\left\{\overline{\Delta}_{\K,f}(h,\cdot),\Delta_{\K,f}(h,\cdot)\right\}$
 and remark that in view of the Lebesgue differentiation theorem
$\Delta^{*}_{\K,f}(h,\cdot)$ and $\overline{\Delta}_{\K,f}(h,\cdot)$ coincide almost everywhere.
Note also, that if $f$ is a continuous function then $\Delta^{*}_{\K,f}(h,\cdot)\equiv \overline{\Delta}_{\K,f}(h,\cdot)$.

%Now we are in a position to introduce the oracle estimator.
Define for $\forall y\in \RR$ the oracle (depending on the underlying function) bandwidth \( h^*_{\K,f}(y) \)
\begin{equation}
\label{def_oracle bandwidth}
 h^*_{\K,f}(y)=
    \sup\big\{h\in [\eps^2, 1]:\quad \sqrt{h}\;\Delta^{*}_{\K,f}(h,y)\leq
      \| \K \|_{\infty}\eps \sqrt{\ln (1/\eps)}\big\}.
\end{equation}
We see that, with the proviso that \( f \in \mathbb{F}_M \), the ``bias'' \( \Delta^{*}_{\K,f}(h,\cdot)\le 2M\|\K\|_1 \), and consequently the set \eqref{def_oracle bandwidth} is not empty for all $\eps\leq \exp{\big\{ -(2M\|\K\|_1\big/\|\K\|_\infty)^2\big\}}$. Here \( \| \K \|_p \), \( 1 \le p \le \infty \), denotes the \( L_p \) norm of \( \K \).
%Some remarks are in order. First, we note that $\Delta^{*}_{\K,f}(h,\cdot)\leq 2M\|\K\|_1$ for any $f\in\mathbb{F}_M$ and, therefore,$ h^*_{\K,f}(y)$ is well-defined for any $y\in \RR$ whenever
%$\eps\leq \exp{\big\{ -(2M\|\K\|_1\big/\|\K\|_\infty)^2\big\}}$, where \( \| \K \|_p \), \( 1 \le p \le \infty \), denotes the \( L_p \) norm of \( \K \).
\begin{comment}
Moreover, Assumption \ref{ass:assumption-on-kernel} ({\it 2}) implies that $\Delta^{*}_{\K,f}(\cdot,y)$ is continuous on $[\eps^2, 1]$, hence
\begin{eqnarray}
\label{eq1:def_oracle bandwidth}
&\text{either}&\quad\sqrt{ h^*_{\K,f}(y)}\;\Delta^{*}_{\K,f}\big( h^*_{\K,f}(y),y\big)=\| \K \|_{\infty}\eps \sqrt{\ln (1/\eps)},
\\*[2mm]
\label{eq2:def_oracle bandwidth}
&\text{or}&
\quad \sqrt{h}\;\Delta^{*}_{\K,f}(h,y)\leq
      \| \K \|_{\infty}\eps \sqrt{\ln (1/\eps)}, \;\;\forall h\in [\eps^2, 1].
\end{eqnarray}
\end{comment}

\begin{comment}

The quantity, similar to the defined above $h^*_{\K,f}(\cdot)$, first appeared in \cite{LepMamSpok97}
in the context of the estimating univariate functions possessing inhomogeneous smoothness.
Some years later this approach has been developed in \cite{LepskiLevit99}, \cite{KLP2001} and  \cite{GoldLep2008}
for the estimation of multivariate function.
In these papers, the interested reader can find a more detailed discussion of the oracle approach. In the present paper we try to adopt the ``ideology'' proposed in the aforementioned papers to the estimation under single index constraint. Our main idea is based on rather simple observation.
\end{comment}

For any $(\theta,h)\in\SS\times[\eps^2, 1]$ define the matrix
$$
E_{(\theta,h)}=\left(
\begin{array}{ll}
h^{-1}\theta_1\quad &h^{-1}\theta_2
\\
-\theta_2\quad &\;\theta_1
\end{array}
\right)
$$
and consider the family of kernel estimators
\begin{equation*}
\label{family of estimators}
\cc F
=
\Big\{ \hat F_{(\ta, h)} (\cdot)
    =
    \det\big(E_{(\ta, h)}\big)\int K\big(E_{(\ta, h)}(t- \cdot)\big) Y_{\eps}(\dd t),\;\; (\theta,h)\in\SS\times[\eps^2, 1]\Big\}.
\end{equation*}
We use the product type kernels $K(u,v)=\K(u)\K(v)$ with a one-dimensional kernel $\K$ obeying Assumption \ref{ass:assumption-on-kernel}. Note also that
$\det\big(E_{(\ta, h)}\big)=h^{-1}$ and
\begin{equation}
\label{eq:distribution-of-kernel-estimator}
\hat F_{(\ta, h)} (\cdot)-\mathbb{E}^{\eps}_F\left[\hat{F}_{(\ta, h)} (\cdot)\right]\quad\sim\quad \mathcal{N}\left(0,\|\K\|^4_2\eps^{2} h^{-1}\right).
\end{equation}
The choice $\theta=\ta^{\circ}$ and $h=h^*:=h^*_{\K,f}(x^{T}\ta^{\circ})$ leads to the ``ideal'' (oracle) estimator $\hat F_{(\ta^{\circ}, h^*)}$, that is the estimator constructed {\it as if} \( \ta^{\circ} \) and \( f \) would be known. Such an ``estimator'' is not available but serves as a quality benchmark, given by the following result.
%First, we note that $\hat F_{(\ta^{\circ}, h^*)} (\cdot)$ is not an estimator in the usual sense, since it depends on the function $F$  to be estimated (more precisely on $(f,\ta^{\circ})$ which determines $F$). The meaning of this estimator is explained by the following result.
\begin{proposition}
\label{prop:risk-of-oracle-estimator}
For any \( (f,\ta^{\circ})\in\mathbb{F}_M\times\SS \), \( \eps\leq \exp{\big\{ -\max [1, (2M\|\K\|_1\big/\|\K\|_\infty)^2]\big\}} \)
 and any \( r\ge 1 \)
\begin{equation*}
    \cc R_{r,x}^{(\eps)} \big(\hat{F}_{(\ta^{\circ}, h^*)}, F\big)\leq \mathfrak{c}_r
\left[\frac{\| \K \|^{4}_{\infty}\;\eps^{2} \ln (1/\eps)}{h^*_{\K,f}(x^{\T}\ta^{\circ})}\right]^{1/2}, \forall x\in[-1/2,1/2]^{2},
\end{equation*}
\end{proposition}
\noindent where $\mathfrak{c}_r=\left[\mathbb{E}\big(1+|\varsigma|\big)^{r}\right]^{1/r},\; \varsigma\sim \mathcal{N}(0,1)$.
 The proof is straightforward and can be omitted.

The meaning of Proposition \ref{prop:risk-of-oracle-estimator} is %can be treated as follows.
that the ``oracle'' knows the exact value of the index $\ta^{\circ}$ and the optimal, up to $\ln(1/\eps)$, bias-variance trade-off $h^*$
between the approximation error caused by $\Delta^{*}_{\K,f}(h^{*},\cdot)$ and the variance, see formula~\eqref{eq:distribution-of-kernel-estimator}, of the kernel estimator from the collection
$\mathcal{F}$.% with bandwidth $h^*$.  %It explains why the ``oracle'' chooses the ``estimator'' $\hat{F}_{(\ta^{\circ}, h^*)}$.

\par Below we will propose an adaptive (not depending of \( \ta^{\circ} \) and \( f \) ) estimator and show that this estimator is as good as the oracle one, i.e. that the risk of that estimator is worse than that of  Proposition \ref{prop:risk-of-oracle-estimator} by a numerical constant only.
\begin{comment}
\par In the next paragraph we propose a ``real'' (based on the observation) estimator $\hat{F}(\cdot)$, which mimics the oracle estimator.
This means that for any $(f,\ta^{\circ})\in\mathbb{F}_M\times\SS$,  $x\in[-1/2,1/2]^{2}$,
$\eps\leq \exp{\big\{ -\max [1, (2M\|\K\|_1\big/\|\K\|_\infty)^2]\big\}}$, and $r>0$
$$
\cc R_{r,x}^{(\eps)} \big(\hat{F}, F\big)\leq C_r\left[\frac{\| \K \|^{4}_{\infty}
\;\eps^{2} \ln (1/\eps)}{h^*_{\K,f}(x^{\T}\ta^{\circ})}\right]^{1/2},
$$
where $C_r$ is an absolute constant independent of the noise level $\eps$ and the underlying function $F$.
The latter result is a local oracle inequality.
The construction of the estimator $\hat{F}(\cdot)$ is based on the data-driven selection from the family $\mathcal{F}$.
\end{comment}

\subsection{Selection rule}\label{sect:SelectionRule} The procedure below is based on a pairwise comparison of the estimators from \( \cc F \) with an auxiliary estimator defined as follows. For any $\ta,\nu\in\SS$  and any $h\in[\eps^2,1]$ introduce the matrices
\begin{equation*}
%\label{matrix for pseudo-estimator}
\overline{E}_{(\theta,h)(\nu,h)}=\left(
\begin{array}{ll}
\frac{(\theta_1+\nu_1)}{2h(1+|\nu^\T\theta|)}\quad &\frac{(\theta_2+\nu_2)}{2h(1+|\nu^\T\theta|)}
\\*[2mm]
-\frac{(\theta_2+\nu_2)}{2(1+|\nu^\T\ta|)}\quad &\;\frac{(\theta_1+\nu_1)}{2(1+|\nu^\T\ta|)}
\end{array}
\right), \quad
%\end{equation*}
%where
%\begin{equation*}
E_{(\theta, h)(\nu,h)}=\left\{
\begin{array}{ll}
\overline{E}_{(\theta, h)(\nu,h)},\;\quad&\nu^\T\ta\geq 0;
\\*[2mm]
\overline{E}_{(-\theta, h)(\nu,h)},\quad&\nu^\T\ta< 0.
\end{array}
\right.
\end{equation*}
It is easy to check that
$
(4h)^{-1}\leq\det\big(E_{(\theta,h)(\nu,h)}\big)\leq (2h)^{-1}.
$
Then, similarly to the construction of the estimators from \( \cc F \) we define a kernel estimator parametrized by \( E_{(\ta, h)(\nu, h)} \)
\begin{equation}\label{pseudo-estimator}
\hat F_{(\ta, h)(\nu, h)}(x)
    =
    \det \big(E_{(\ta, h)(\nu,h)}\big)
        \int K( E_{(\ta, h)(\nu,h)}(t-x))Y_{\eps}(\dd t).
\end{equation}

Put $\Lambda(\K,Q)=8\sqrt{\ln{(1+ 2Q\|\K\|_\infty)}}+50$ and let for any $\eta\in (0,1]$
$$
\TH(\eta)=2\| \K \|^{2}_{\infty}\left[\Lambda(\K,Q)+\sqrt{4r+2}+1\right]\eps\sqrt{\eta^{-1}\ln (1/\eps)}.
$$

Set
$\mathcal{H}_\eps=\big\{h_k=2^{-k},\; k = 0,1, \ldots\big\}\cap[\eps^2,1]$ and define for any \(\ta\in \SS \)
and  \(h \in \mathcal{H}_\eps \)
\begin{equation}\label{eq:proced_first-step}
R_{(\ta, h)}(x)
=
\sup_{\eta\in\mathcal{H}_\eps:\;  \eta\le h}
\Big\{
    \sup_{\nu\in\SS}
        \big|\hat F_{(\ta, \eta)(\nu, \eta)}(x) - \hat F_{(\nu, \eta)}(x) \big|   -\TH(\eta)
    \Big\}.
\end{equation}

\noindent For any $x\in[-1/2,1/2]^{2}$ introduce the random set
\begin{eqnarray*}
\cc P (x)
=
    \big\{   (\ta, h)\in\SS\times\mathcal{H}_\eps : R_{(\ta, h)}(x) \le 0 \big\},
\end{eqnarray*}
and let  $\tilde{h}=
\max\big\{h:\;\; (\theta,h)\in\cc P (x)\big\}$ if  $\cc P (x)\neq\emptyset.$  Note that there exists $\vartheta\in\SS$
such that $(\vartheta,\tilde{h})\in\mathcal{P}(x)$, since the set $\mathcal{H}_\eps$ is finite.
Define
\begin{eqnarray*}
\hat{\ta}=\left\{
\begin{array}{ll}
(1,0)^{\T},\quad & \cc P (x)=\emptyset;
\\*[1mm]
\theta \text{   s.t. } (\theta,\tilde{h})\in\cc P (x),\quad& \cc P (x)\neq\emptyset.
\end{array}
\right.
\end{eqnarray*}
If $\hat{\theta}$ is not unique, let us make any measurable choice. In particular, if $\hat{\Theta}:=\big\{\theta\in\SS:\;(\theta,\tilde{h})\in\cc P (x)\big\}$ one can choose $\hat{\theta}$ as a vector belonging to
$\hat{\Theta}$ with the smallest first coordinate. The measurability of this choice follows from the fact that the mapping $\ta\mapsto R_{(\ta, h)}(x)$ is almost surely continuous  on $\SS$. This continuity, in its turn, follows from Assumption \ref{ass:assumption-on-kernel} ({\it 2}), bound (\ref{eq100:proof-lemma-gauss-on-matrices}) for Dudley's entropy integral proved in Lemma \ref{lem:gauss-on-matrices} below and the condition $f\in\mathbb{F}_M$.
Define
\begin{equation}\label{eq:proced_second-step}
\hat{h}=\sup\left\{h\in\mathcal{H}_\eps:\;\;\left|\hat F_{(\hat{\ta},h)}(x)- \hat F_{(\hat{\ta}, \eta)}(x)\right|
                \leq\TH(\eta),\;\;\forall \eta\leq h,\;\eta\in\mathcal{H}_\eps\right\}
\end{equation}
and put as a final  estimator
\( \hat F (x)=\hat F_{(\hat \ta, \hat h)}(x)\).

\par The proposed above procedure belongs to the stream of pointwise adaptive procedures originating from \cite{Lep1990}. Indeed, the second step determined by \eqref{eq:proced_second-step} for the ``frozen'' \( \hat \ta \) is exactly the procedure of \cite{Lep1990} which was originally developed in the framework of the univariate GWN model. There is a rather vast literature on that topic, we mention \cite{BauerHohageMunk} adapted the method of  \cite{Lep1990} for the choice of the parameter for iterated Tikhonov regularization in nonlinear inverse problems, \cite{BertinRivoirard} showed the maxiset optimality of that procedure for bandwidth selection under the \( sup \) norm losses, \cite{Chichignoud} used it for selecting among local bayesian estimators, \cite{Gaiffas} studied the problem of pointwise estimation in random design Gaussian regression, \cite{Serdyukova} investigated a heteroscedastic Gaussian regression under noise misspecification, among many others.

\par The application of \cite{Lep1990} requires some sort of ordering on the set of estimators, for instance in~\eqref{eq:proced_second-step} as soon as \( \hat \ta \) is fixed it is due to the monotonicity of the ``bias'' \( \Delta^{*}_{\K,f}(\cdot,y) \). However, when the projection direction is unknown no natural order on \( \cc F \) is available. This problem is similar to the one arising in generalizations of the pointwise adaptive method for multivariate (anisotropic) settings, see for developments in that direction \cite{LepskiLevit99}, \cite{KLP2001} and  \cite{GoldLep2009}. Usually the aforementioned issue requires to introduce an auxiliary estimator and construct a procedure carefully capturing the ``incomparability'' of the estimators. In the considered set-up it is realized by the first step of procedure with \( R_{(\ta, h)}(x) \) given by \eqref{eq:proced_first-step}.

\subsection{Oracle inequalities}\label{subsec:OI}

Throughout the paper   we assume that

$$
\eps\leq \exp{\big\{ -\max [1, (2M\|\K\|_1\big/\|\K\|_\infty)^2]\big\}}.
$$

\begin{theorem}
\label{th:local-oracle-inequality}
For any $(f,\ta^{\circ})\in\mathbb{F}_M\times\SS$,  $x\in[-1/2,1/2]^{2}$ and any $r\ge 1$
$$
\cc R_{r,x}^{(\eps)} \Big(\hat{F}_{(\hat{\theta}, \hat{h})}, F\Big)
    \leq
    C_{r,1}(Q,\K)\sqrt{\frac{\|\K\|^4_\infty\eps^{2} \ln (1/\eps)}{h^*_{\K,f}(x^{T}\ta^{\circ})}}
        +
            C_{r,2}(M,Q,\K)\| \K \|^2_{\infty}\eps\sqrt{\ln (1/\eps)}.
$$
The constants $C_{r,1}(Q,\K)$ and $C_{r,2}(M,Q,\K)$ are given in the beginning of the proof.

\end{theorem}

As already mentioned, the global oracle inequality is
obtained by integrating the local oracle inequality. For ease of notation, we write \( r(\eps) = C_{r,2}(M,Q,\K)\| \K \|^2_{\infty}\,\eps\sqrt{\ln (1/\eps)} \) and $C_r=C_{r,1}(Q,\K)$. It follows from Jensen's inequality and Fubini's theorem that
$$
\cc R_r^{(\eps)} (\hat F, F)
\leq
%\expect{F} \big\| \hat F (\cdot) - F(\cdot) \big\|^{r}_r
\big\|\cc R_{r,\cdot}^{(\eps)} (\hat F, F)\big\|_r
\leq
C_r \left\{ \int_{[-1/2,1/2]^{2}}\left[\frac{\|\K\|^4_\infty\eps^{2}
 \ln (1/\eps)}{h^*_{\K,f}(x^{T}\ta^{\circ})}\right]^{\frac{r}{2}}\dd x \right \}^{\frac{1}{r}} + r(\eps) .
$$
Integration by substitution gives:
$$
\int_{[-1/2,1/2]^{2}}\left[\frac{\|\K\|^4_\infty\eps^{2}
 \ln (1/\eps)}{h^*_{\K,f}(x^{T}\ta^{\circ})}\right]^{\frac{r}{2}}\dd x\leq \int_{-1/2}^{1/2}\left[\frac{\|\K\|^4_\infty\eps^{2}
 \ln (1/\eps)}{h^*_{\K,f}(z)}\right]^{\frac{r}{2}}\dd z
$$
leading to the following result.
\begin{theorem}
\label{th:global-oracle-inequality}
For any $(f,\ta^{\circ})\in\mathbb{F}_M\times\SS$ and any $r\ge1$
$$
\cc R_r^{(\eps)} \Big(\hat{F}_{(\hat{\theta}, \hat{h})}, F\Big)\leq C_{r,1}(Q,\K)\left\|\sqrt{\frac{\|\K\|^4_\infty\eps^{2}
 \ln (1/\eps)}{h^*_{\K,f}(\cdot)}}\right\|_r+C_{r,2}(M,Q,\K)\| \K \|^2_{\infty}\eps\sqrt{\ln (1/\eps)}.
$$

\end{theorem}

\section{Adaptation}
\label{subsec:adaptive-estimation}

In this section with the use of the local oracle inequality from Theorem~\ref{th:local-oracle-inequality} we solve the problem of pointwise adaptive estimation over a collection of H\"older classes. Then, we turn to the problem of adaptive estimating the entire function over a collection of Nikol'skii classes with the accuracy of an estimator measured under the~$L_r$~risk. That is done with the help of the global oracle inequality given in Theorem~\ref{th:global-oracle-inequality}.

 \par Throughout this section we will assume that the kernel $\K$ satisfies additionally Assumption \ref{ass2:assumption-on-kernel} below.
Introduce the following notation: for any $a>0$ let $m_a \in\mathbb{N}$ be the maximal integer
strictly less than $a$.

\begin{assumption}
\label{ass2:assumption-on-kernel}
There exists $\bob>0$ such that
$$
\int z^{j}K(z)\dd z=0,\;\;\forall j=1,\ldots,m_{\bob}.
$$

\end{assumption}

\subsection{Pointwise adaptation}\label{subsec:PointwiseAdaptation}

%We start this section with the definition of  the H\"older class of functions.
Let us firstly recall the definition of H\"olderian functions.
\begin{definition}
\label{def:holder-class}

Let \( \beta>0 \) and \( L>0 \). A function $g:\mathbb{R}\to\RR$ belongs to the H\"older class $\mathbb{H}(\beta,L)$ if $g$ is $m_\beta$-times continuously
differentiable, $\|g^{(m)}\|_\infty\leq L,\;\forall m\leq m_\beta,$  and
$$
\left|g^{(m_\beta)}(t+h)-g^{(m_\beta)}(t)\right|\le L h^{\beta-m_\beta},\;\;\forall t \in \RR \; \text{and} \; h>0.
$$

\end{definition}

\par The aim is to estimate the function $F(x)$ at a given point $x\in [-1/2,1/2]^2$ under the additional assumption that $F\in\mathbb{F}(\bob):=\bigcup_{\beta\leq\bob}\bigcup_{L>0}\mathbb{F}_2(\beta,L)$, where
\begin{equation}\label{Def:SI-class_local}
    \mathbb{F}_d(\beta,L)=\left\{F:\RR^d\to\RR\;|\; F(z)=f\big(z^{\T}\theta\big),\;f\in\mathbb{H}(\beta,L),\;\theta\in\SSd\right\},
\end{equation}
$d\geq 2$ is the dimension and $\bob$ is the constant  from Assumption \ref{ass2:assumption-on-kernel}, which can be arbitrary but must be chosen {\it a priory}.
% We will see that $\bob$ can be an arbitrary number but it must be chosen {\it a priory}.

\begin{theorem}
\label{th:pointwise-adaptation}
Let $\bob>0$ be fixed and let Assumptions \ref{ass:assumption-on-kernel} and \ref{ass2:assumption-on-kernel} hold. Then, for any
$\beta\leq\bob$, $L>0$ and $x\in[-1/2,1/2]^{2}$, we have
$$
\sup_{F\in\mathbb{F}_2(\beta,L)}\cc R_{r,x}^{(\eps)} \Big(\hat{F}_{(\hat{\theta}, \hat{h})}, F\Big)
\leq  \|\K \|^2_{\infty}\left[%2^{\frac{1}{2\beta+1}}
C_{r,1}(Q,\K)\psi_\eps(\beta,L)
+C_{r,2}(L,Q,\K)\, \eps\sqrt{\ln (1/\eps)}\right],
$$
where $\psi_\eps(\beta,L)=L^{\frac{1}{2\beta+1} }\left(\eps\sqrt{\ln(1/\eps)}\right)^{\frac{2\beta}{2\beta+1}}$.

\smallskip

\par Moreover, for any $\beta,L>0$, \( r \ge 1 \), \( x \in [-1/2, 1/2]^d \) with \( d\ge 2 \)  and  any $\eps>0$ small enough,
$$
\inf_{\tilde{F}}\sup_{F\in\mathbb{F}_d(\beta,L)}\cc R_{r,x}^{(\eps)} \Big(\tilde{F}, F\Big)
\geq \kappa \psi_\eps(\beta,L),
$$
where infimum is over all possible estimators. Here $\kappa$ is a numerical constant  independent of $\eps$ and $L$.
\end{theorem}
We conclude that the estimator $\hat{F}_{(\hat{\theta}, \hat{h})}$ is minimax adaptive with respect to the collection of classes
$\big\{\mathbb{F}_d(\beta,L),\;\; \beta\le \bob,\; L>0\big\}$. As already mentioned, this result is quite surprising. Indeed, if for example,
the directional vector $\theta=(1,0)^{\T}$, i.e. is known, then $\mathbb{F}(\beta,L)=\mathbb{H}(\beta,L)$ and the considered estimation problem
 can be easily reduced to estimation of $f$ at a given point in the univariate Gaussian white noise model. As it is shown in \cite{Lep1990}
the adaptive estimator over the collection  $\big\{\mathbb{H}(\beta,L),\;\; \beta\leq \bob,\; L>0\big\}$ does not exist.

Also, we would like to emphasize that the lower bound result given by the second assertion of the theorem is proved for arbitrary dimension. As to the proof of the first statement of the theorem it is based on the evaluation of the uniform, over $\mathbb{H}_d(\beta,L)$, lower bound for
$
h^*_{\K,f}(\cdot)
$
and on the application of Theorem \ref{th:local-oracle-inequality}. We note also that the upper bound for the minimax risk given in Theorem \ref{th:pointwise-adaptation} was earlier given in \cite{GoldLep2008}, but the estimation procedure used there is completely
different from our selection rule.

\subsection{Adaptive estimation under the \( L_r \) losses}\label{sect:Global_Adapt}

We start this section with the definition of  the Nikol'skii class of functions.

\begin{definition}
\label{def:nikolskii-class}

Let \( \beta>0 \), \( L>0 \) and $p\in [1,\infty)$ be fixed. A function $g:\mathbb{R}\to\RR$ belongs to the Nikol'skii class $\mathbb{N}_p(\beta,L)$,
if $g$ is $m_\beta$-times continuously
differentiable and
\begin{eqnarray*}
&&\left( \int_{\RR}\left|g^{(m)}(t)\right|^p \dd t \right)^{\frac{1}{p} }\le L,\quad \forall m=1,\ldots, m_\beta;
\\
&&\left( \int_{\RR}\left|g^{(m_\beta)}(t+h)-g^{(m_\beta)}(t)\right|^p \dd z \right)^{\frac{1}{p} }
\leq Lh^{\beta-m_\beta},\;\;\forall h>0.
\end{eqnarray*}

\end{definition}
\noindent Later on  we assume that $\mathbb{N}_p(\beta,L)=\mathbb{H}(\beta,L)$ if $p=\infty$.

\par Here the target of estimation is the entire function $F(\cdot)$  under the assumption that
$ F \in\mathbb{F}_p(\bob):=\bigcup_{\beta\leq\bob}\bigcup_{L>0}\mathbb{F}_{2,p}(\beta,L)$, where
\begin{equation}\label{Def:SI-class_global}
    \mathbb{F}_{d,p}(\beta,L)=\left\{F:\RR^d\to\RR\;|\; F(z)=f\big(z^{\T}\theta\big),\;f\in\mathbb{N}_p(\beta,L),\;\theta\in\SSd\right\}.
\end{equation}

\begin{comment}
Let us briefly discuss the applicability of  Theorem \ref{th:global-oracle-inequality} which requires
that $f\in\mathbb{F}_M$. In order to guarantee it we will assume that
 $\beta p> 1$. The latter assumption is standard in estimation of functions possessing inhomogeneous smoothness, see for example,
\cite{DJKP}, \cite{LepMamSpok97}, \cite{KLP2008}.
If  $\beta p> 1$ the embedding
$\mathbb{N}_p(\beta,L)\subset\mathbb{H}(\beta-1/p,cL)$ with some absolute constant $c>0$ guarantees that $f\in\mathbb{F}_{M},\; M=cL$
and Theorem \ref{th:global-oracle-inequality} is applicable.
\end{comment}

\begin{theorem}
\label{th:global-adaptation}
Let $\bob>0$ be fixed and let Assumptions \ref{ass:assumption-on-kernel} and \ref{ass2:assumption-on-kernel} hold. Then, for any
 $L>0$, $p>1$, $p^{-1}<\beta\leq\bob$ and  $r\geq 1$,
$$
\sup_{F\in\mathbb{F}_{2,p}(\beta,L)}\cc R_{r}^{(\eps)} \Big(\hat{F}_{(\hat{\theta}, \hat{h})}, F\Big)
\leq
\|\cc K\|^{2}_\infty \Big[ \kappa C_{r,1}(Q,\K)\varphi_\eps(\beta,L,p)+C_{r,2}(L,Q,\K)\eps\sqrt{\ln (1/\eps)}\Big],
$$
where
$$
\varphi_\eps(\beta,L,p)=\left\{
\begin{array}{lll}
L^{\frac{1}{2\beta+1}}\left(\eps\sqrt{\ln(1/\eps)}\right)^{\frac{2\beta}{2\beta+1}},\quad& (2\beta+1)p>r;
\\
L^{\frac{1}{2\beta+1}}\left(\eps\sqrt{\ln(1/\eps)}\right)^{\frac{2\beta}{2\beta+1}}\big[\ln(1/\eps)\big]^{\frac{1}{r}},\quad& (2\beta+1)p=r;
\\
L^{\frac{1/2-1/r}{\beta-1/p+1/2}}
\left(\eps\sqrt{\ln (1/\eps)}\right)^{\frac{\beta-1/p+1/r}{\beta-1/p+1/2}},\quad& (2\beta+1)p< r.
\end{array}
\right.
$$
The constant $\kappa$ is independent of $\eps$, $L$ and $\K$.

\end{theorem}

Let us make some remarks. First, note that $\mathbb{F}_{2,p}(\beta,L)\supset\mathbb{N}_{p}(\beta,L)$.
Indeed, the class  $\mathbb{N}_{p}(\beta,L)$ can be viewed as the class of functions $F$ satisfying $F(\cdot)=f(\ta^\T\cdot)$ with
$\ta=(1,0)^{\T}$. Then, the problem of estimating such (2-variate) functions  can be reduced  to the estimation
of univariate functions observed
in the one-dimensional GWN model. In view of this remark the rate of convergence for the latter problem
(which  can be found for example in \cite{DJKP,DelyonJuditski1996} ) is the lower bound for the minimax risk  defined on
$\mathbb{F}_{2,p}(\beta,L)$. Under assumption $\beta p>1$ this rate of convergence is given by
$$
\phi_\eps(\beta,L,p)=\left\{
\begin{array}{lll}
L^{\frac{1}{2\beta+1}}\eps^{\frac{2\beta}{2\beta+1}},\quad& (2\beta+1)p> r;
\\
L^{\frac{1}{2\beta+1}}\left(\eps\sqrt{\ln(1/\eps)}\right)^{\frac{2\beta}{2\beta+1}},\quad& (2\beta+1)p=r;
\\
L^{\frac{1/2-1/r}{\beta-1/p+1/2}}
\left(\eps\sqrt{\ln (1/\eps)}\right)^{\frac{\beta-1/p+1/r}{\beta-1/p+1/2}},\quad& (2\beta+1)p< r.
\end{array}
\right.
$$
The minimax rate of convergence in the case $(2\beta+1)p=r$ remains an open problem, and the rate presented in the middle line  above
is only the lower asymptotic bound for the minimax risk. Therefore the proposed estimator $\hat{F}_{(\hat{\theta}, \hat{h})}$
is adaptive  whenever  $(2\beta+1)p<r$.
In the case $(2\beta+1)p\geq r$ we loose only a logarithmic factor with respect to the optimal rate
and, as mentioned in Introduction, the construction of  adaptive estimator over a collection
$\big\{\mathbb{F}_{2,p}(\beta,L), \;\beta>0,\;L>0\big\}$
in this case remains an open problem.

\section{Proofs}
\label{sec:proofs}

\subsection{Proof of Theorem \ref{th:local-oracle-inequality}}

   The section starts with % a list of
   the constants used in the statement of the theorem as well as technical lemmas whose proofs are postponed to Appendix.

\paragraph{Constants}
\begin{eqnarray*}
&&C_{r,1}(Q,\K)
    =
    8  \left[\Lambda(\K,Q)+\sqrt{4r+2}+1\right]
            +
        \mathfrak{c}_r  \left[ (2 + \sqrt{2}) \Lambda(\K,Q) + 2 \right] +1;
\\*[2mm]
&&C_{r,2}(M,Q,\K)=2^{1/r}\left[2M+\Lambda(\K,Q) \mathfrak{c}_{2r} \right] .
\end{eqnarray*}

\subsubsection{Auxiliary results}
For any $\theta,\nu\in\SS$ and $h\in[\eps^2,1]$ denote
\begin{eqnarray*}
%\label{pseudo-smoother}
S_{(\ta, h)(\nu, h)}(x)
    &=&
    \det \big(E_{(\ta, h)(\nu,h)}\big)
        \int K(E_{(\ta, h)(\nu,h)}(t-x))F(t)\dd t ,
\\
S_{(\ta, h)}(x)
    &=&
    \det \big(E_{(\ta, h)}\big)
        \int K(E_{(\ta, h)}(t-x))F(t)\dd t .
\end{eqnarray*}
For ease of notation, we write $h^*_f=h^*_{\K,f}(x^{\T}\ta^{\circ})$.

\begin{lemma}
\label{lem:bounds for bias} Grant Assumption \ref{ass:assumption-on-kernel}. Then, for any $\nu\in\SS$ and any $\eta,h\in [\eps^2,1]$
satisfying $\eta\leq h\leq 2^{-1}h^*_{f}$, one has
\begin{eqnarray*}
&&\left|S_{(\ta^{\circ}, h)(\nu, h)}(x)-S_{(\nu, h)}(x)\right|\leq 2(h^*_f)^{-1/2}\| \K \|^{2}_{\infty}\;\eps \sqrt{\ln (1/\eps)};
\\[2mm]
&&\left|S_{(\nu, h)}(x)-S_{(\nu, \eta)}(x)\right|\leq 2(h^*_f)^{-1/2}\| \K \|^{2}_{\infty}\;\eps \sqrt{\ln (1/\eps)};
\\[2mm]
&&\left|S_{(\ta^{\circ}, h)} - F(x)\right|\leq (h^*_f)^{-1/2}\| \K \|_{\infty}\;\eps \sqrt{\ln (1/\eps)}.
\end{eqnarray*}

\end{lemma}

\smallskip

Let $\mathcal{E}_{a,A},\; 0<a,A<\infty,$ be a set of $2\times2$ matrices such that
$$
\left|\det(E)\right|\geq a,\quad |E|_\infty\leq A,\;\;\forall E\in\mathcal{E}_{a,A}.
$$
Here $|E|_\infty = \max_{i,j} |E_{i,j}|$ denotes the supremum norm, the maximum absolute value entry of the matrix \( E \). Later on without loss of generality we will assume that $a\leq A,\; A\geq 1$.

Assume that the function $ \cc L:\mathbb{R}^{2}\to\mathbb{R}$ is compactly supported on $[-1/2,1/2]^{2}$,  $\int \cc L=1$ and satisfies the Lipschitz condition
$$
\left| \cc L(u)- \cc L(v)\right|\leq \Upsilon |u-v|_2,\;\;\forall u,v\in \RR^2,
$$
where $|\cdot|_2$ is the Euclidian norm. Let $y\in \RR^2$ be fixed. On the parameter set  $\mathcal{E}_{a,A}$ let a Gaussian random function be defined by
$$
\zeta_{y}(E)= \| \cc L\|^{-1}_2\sqrt{\left|\det(E)\right|}\int  \cc L\big(E(u-y)\big)W(\dd u).
$$

Put \( \mathbf{c}(a,A) =4\sqrt{2}  \left[ \ln(A\vee\{A/a\}^2) + 2 \ln{(1+ \sqrt{2} \Upsilon)} \right]^{1/2} + 29 \)
%$\mathbf{c}(a,A)=8 \left[ \ln(A\vee\{A/a\}) + \ln{(4+2\sqrt{2}\Upsilon)} \right]^{1/2} + 20 $
and  \( \mathfrak{c}_q=\left(\mathbb{E}\big(1+|\varsigma|\big)^{q}\right)^{1/q} \), where  \( \varsigma\;\sim\;\norm{0}{1} \) .
\begin{lemma}
\label{lem:gauss-on-matrices}
For any \( z>0 \)
\begin{equation*}
\mathbb{P}\Big\{\sup_{E\in\mathcal{E}_{a,A}}\left|\zeta_{y}(E)\Big| \ge \mathbf{c}(a,A)+z\right\}\le \PP\{ |\varsigma| \ge z \}\le e^{-\frac{z^2}{2}}.
\end{equation*}
Moreover, for any \( q\ge 1 \)
\begin{equation*}
    \left(\mathbb{E}\Big[\sup_{E\in\mathcal{E}_{a,A}}\big|\zeta_{y}(E)\big|\Big]^{q}\right)^{1/q}
    \le
    \mathfrak{c}_q \mathbf{c}(a,A)
%\left(1+\mathfrak{c}_q\right)\mathbf{c}(a,A).
%\left(1+(2\pi)^{1/(2q)}\mathfrak{c}_q\right)\mathbf{c}(a,A).
\end{equation*}

\end{lemma}

\subsubsection{Proof of Theorem \ref{th:local-oracle-inequality}}

Let $h^*\in\mathcal{H_\eps}$ be such that
$
h^*\leq 2^{-1}h^*_{f}<2h^*
$.
Introduce the random events
$$
\cc A=\left\{(\ta^{\circ},h^*)\in\mathcal{P}(x)\right\},\quad
\cc B=\left\{\hat{h}\ge h^*\right\},\quad \cc C=\cc A \cap \cc B,
$$
and let \( \overline{\cc C} \) denote the event complimentary to $\mathcal{C}$. We split the proof into two steps.

\paragraph{Risk computation under \( \cc C \)} The triangle inequality gives
\begin{eqnarray}
\label{eq0:proof-of-theorem-local}
\left|\hat F_{(\hat{\ta},\hat{h})}(x)-F(x)\right|&\leq&\left|\hat F_{(\hat{\ta},\hat{h})}(x)- \hat F_{(\hat{\ta},h^*)}(x)\right|+
\left|\hat F_{(\ta^{\circ},h^*)(\hat{\ta},h^*)}(x)-\hat F_{(\hat{\ta},h^*)}(x) \right|
\nonumber\\*[2mm]
&\;&+\left|\hat F_{(\ta^{\circ},h^*)(\hat{\ta},h^*)}(x)-\hat F_{(\ta^{\circ},h^*)}(x) \right|+\left|\hat F_{(\ta^{\circ},h^*)}(x)-F(x)\right|.
\end{eqnarray}

$1^0.\;$ Since $h^*\geq 4^{-1}h^*_f$ the definition of $\hat{h}$ yields
\begin{eqnarray}
\label{eq1:proof-of-theorem-local}
\left|\hat F_{(\hat{\ta},\widehat{h})}(x)- \hat F_{(\hat{\ta},h^*)}(x)\right|\mathrm{1}_{\mathcal{B}}\leq \TH(h^*)\leq\TH\big(h_f^*/4\big).
\end{eqnarray}

Let us make some remarks. Note that $E_{(\theta, h)(\nu,h)}=\pm E_{(\nu, h)(\theta,h)}$ for any $\theta,\nu$ and $h$. Hence, we conclude that
$\hat F_{(\ta^{\circ},h^*)(\hat{\ta},h^*)}(\cdot)\equiv \hat F_{(\hat{\ta},h^*)(\ta^{\circ},h^*)}(\cdot)$ since $\mathcal{K}$ is symmetric, see Assumption \ref{ass:assumption-on-kernel}.
Next, we note that obviously  $\mathcal{A}\subseteq\{\mathcal{P}(x)\neq\emptyset\}$
 and, moreover, $\mathcal{A}\subseteq\big\{\tilde{h}\geq h^*\big\}$ in view of the definition
of $\tilde{h}$. Lastly, $\big(\hat{\theta},\tilde{h}\big)\in\cc P(x)$ by definition that means $R_{(\hat{\ta}, \tilde{h})}(x)\leq 0$. Consequently,
\begin{eqnarray}
\label{eq2:proof-of-theorem-local}
\hskip-1cm \left|\hat F_{(\ta^{\circ},h^*)(\hat{\ta},h^*)}(x)-
\hat F_{(\ta^{\circ},h^*)}(x) \right|\mathrm{1}_{\mathcal{A}}
&=&\left|\hat F_{(\hat{\ta},h^*)(\ta^{\circ},h^*)}(x)-\hat F_{(\ta^{\circ},h^*)}(x) \right|
\mathrm{1}_{\mathcal{A}}
\nonumber\\*[2mm]
&\leq&\TH(h^*) \leq \TH\big(h_f^*/4\big).
\end{eqnarray}

$2^0.\;$Introduce the following notations. For any $\ta,\nu\in\SS$ and $h\in[\eps^2,1]$ set
\begin{eqnarray*}
\xi _{(\ta, h)(\nu,h)}(x)
    &=&
    \|K\|^{-1}_2\sqrt{\det\big( E_{(\ta, h)(\nu,h)}\big) }
        \int K\big( E_{(\ta,h)(\nu,h)}(t-x)\big)W (\dd t);
\\
\xi _{(\ta, h)}(x)
    &=&
    \|K\|^{-1}_2\sqrt{\det\big( E_{(\ta, h)}\big)}
        \int K\big( E_{(\ta,h)}(t-x)\big)W (\dd t).
\end{eqnarray*}
We remark that  $\big| E_{(\ta, h)}\big|_\infty\leq h^{-1}$ and $\big|E_{(\theta,h)(\nu,h)}\big|_\infty\leq h^{-1}$. Moreover,
\begin{eqnarray}
\label{eq3:proof-of-theorem-local}
&&(4h)^{-1}\leq\det\big(E_{(\theta,h)(\nu,h)}\big)\leq (2h)^{-1},\qquad\det\big( E_{(\ta, h)}\big)=h^{-1}.
\end{eqnarray}
 Since  $h\in[\eps^2,1]$, we assert that
\begin{equation}
\label{eq4:proof-of-theorem-local}
E_{(\theta,h)(\nu,h)}, E_{(\ta, h)}\in \mathcal{E}_{\frac{1}{4},\frac{1}{\eps^{2}}},\quad \forall \ta,\nu\in\SS,\;\forall h\in[\eps^2,1].
\end{equation}
We note also that
for any $\ta,\nu\in\SS$ and $h\in[\eps^2,1]$
\begin{eqnarray*}
&&\left|\hat F_{(\ta^{\circ},h^*)(\hat{\ta},h^*)}(x)-\hat F_{(\hat{\ta},h^*)}(x)\right|
%\\*[2mm]&\le&
\leq\left|S_{(\ta^{\circ},h^*)(\hat{\ta},h^*)}(x)-S_{(\hat{\ta},h^*)}(x)\right|
\\*[2mm]
&&+
\eps\|K\|_2\sqrt{\det\big( E_{(\ta^{\circ},h^*)(\hat{\ta},h^*)}\big) }\left|\xi _{(\ta^{\circ},h^*)(\hat{\ta},h^*)}(x)\right|
+\eps\|K\|_2\sqrt{\det\big( E_{(\hat{\ta},h^*)}\big) }\left|\xi _{(\hat{\ta},h^*)}(x)\right|.
\end{eqnarray*}
We obtain from the first assertion of Lemma  \ref{lem:bounds for bias} with $\nu=\hat{\ta},\;h=h^*$,
(\ref{eq3:proof-of-theorem-local}) and (\ref{eq4:proof-of-theorem-local})
\begin{eqnarray}
\label{eq5:proof-of-theorem-local}
\left|\hat F_{(\ta^{\circ},h^*)(\hat{\ta},h^*)}(x)-\hat F_{(\hat{\ta},h^*)}(x)\right|
&\leq&
\frac{2\,\| \K \|^{2}_{\infty}}{\sqrt{h^*_f}} \,\eps \sqrt{\ln (1/\eps)}
    +
     \frac{2 + \sqrt{2}}{\sqrt{h^*_f}}\;\| \K \|^{2}_{2}\,\eps \sqrt{\ln (1/\eps)}\;\zeta_\eps(x)
\nonumber\\
&\leq&
\frac{\| \K \|^{2}_{\infty}}{\sqrt{h^*_f}} \, \eps \sqrt{\ln (1/\eps)}\,\big[2+(2+\sqrt{2})\zeta_\eps(x)\big],
\end{eqnarray}
where we denoted
$$
\zeta_\eps=\left[\ln{(1/\eps)}\right]^{-1/2}
\displaystyle{\sup_{E\in \mathcal{E}_{\frac{1}{4},\frac{1}{\eps^{2}}}}}\left|\zeta_x\big(E\big)\right|.
$$
We have also used that $2h^*\leq h^*_f<4h^*$.

\smallskip

$3^0.\;$ We get in view of the third assertion of Lemma  \ref{lem:bounds for bias} that
\begin{eqnarray}
\label{eq6:proof-of-theorem-local}
\left|\hat F_{(\ta^{\circ},h^*)}(x)-F(x)\right|
&\leq&  \sqrt{1/h_f^*}\;\| \K \|_{\infty}\;\eps
\sqrt{\ln (1/\eps)}+\sqrt{1/h^*}\;\| \K \|^{2}_{2}\;\eps|\varsigma|
\nonumber\\*[2mm]
&\leq&\sqrt{1/h_f^*}\;\| \K \|^2_{\infty}\;\eps \sqrt{\ln (1/\eps)}\big(1+2|\varsigma|\big),
\end{eqnarray}
where $\varsigma\;\sim\;\mathcal{N}(0,1).$

\smallskip

$4^0.\;$ We obtain from (\ref{eq0:proof-of-theorem-local}), (\ref{eq1:proof-of-theorem-local}), (\ref{eq2:proof-of-theorem-local}), (\ref{eq5:proof-of-theorem-local}) and (\ref{eq6:proof-of-theorem-local}) and the second assertion of Lemma \ref{lem:gauss-on-matrices}
with $ \cc L=K,\;a=1/4,\; A=\eps^{-2}$ and  $q=r$, noting that $\Upsilon=\sqrt{2}Q\|\K\|_\infty$ ,
\begin{eqnarray}
\label{eq7:proof-of-theorem-local}
\left\{\mathbb{E}\left|\hat F_{(\hat{\ta},\hat{h})}(x)-F(x)\right|^r\mathrm{1}_{\mathcal{C}}\right\}^{1/r}
%\nonumber\\*[2mm]
&\le&
 2\TH\big(h^*_f/4\big)
 +
  \left[ (2 + \sqrt{2}) \Lambda(\K,Q) \mathfrak{c}_r + 2 \mathfrak{c}_r +1 \right]
%\left[
 %   \Lambda(\K,Q)  (2 + \sqrt{2}) \left(1+(2\pi)^{1/(2r)}\mathfrak{c}_r\right)
    %        +    2\mathfrak{c}_r+1 \right]
\frac{\| \K \|^{2}_{\infty}}{\sqrt{h^*_f}}\eps \sqrt{\ln (1/\eps)}
\nonumber\\
&\le&
 C_{r,1} \frac{\| \K \|^{2}_{\infty}}{\sqrt{h^*_f}} \eps \sqrt{\ln (1/\eps)},
\end{eqnarray}
Here we have also used that
$$
\sup_{\eps\leq e^{-1}}
    \left[
    \frac{4\sqrt{2}  \left\{ 2 \sqrt{\ln{(2/\eps)}}+ \sqrt{2\ln(1+2Q\|\K\|_\infty)} \right\} + 29}{\sqrt{\ln{(1/\eps)}}}
    \right]\leq \Lambda(\K,Q).
$$

\paragraph{Risk computation under \( \overline{\cc C} \)}

 Since $f\in\mathbb{F}_M$ one can easily evaluate the discrepancy between the adaptive estimator and the value of function
$$
\left|\hat F_{(\hat{\ta},\hat{h})}(x)-F(x)\right|\leq M\big(1+\|K\|_1\big)+\eps\|K\|_2\sqrt{\det\big( E_{(\hat{\ta},\hat{h})}\big) }\left|\xi _{(\hat{\ta},\hat{h})}(x)\right|.
$$
We obtain in view of (\ref{eq3:proof-of-theorem-local}) and (\ref{eq4:proof-of-theorem-local}), taking into account that $\hat{h}>\eps^{2}$,
$$
\left|\hat F_{(\hat{\ta},\hat{h})}(x)-F(x)\right|\leq M\big(1+\|\K\|^{2}_1\big)+\|\K\|^2_2\sqrt{\ln (1/\eps)}\zeta_\eps.
$$
Thus, applying the second assertion of Lemma \ref{lem:gauss-on-matrices}
with $ \cc L=K,\;a=1/4,\; A=\eps^{-2}$,  $\Upsilon=\sqrt{2}Q\|\K\|_\infty$ and $q=2r$, we get
\begin{equation*}
    \left[\mathbb{E}_F^{\eps}\left|\hat F_{(\hat{\ta},\hat{h})}(x)-F(x)\right|^{2r}\right]^{1/2r}
    \le
    \left[2M+\Lambda(\K,Q) \mathfrak{c}_{2r} \right] \|\K\|^2_\infty  \sqrt{\ln (1/\eps)}.
\end{equation*}
Here it is used that $1\le\|\K\|_1\le \|\K\|_2\le \|\K\|_\infty$ due to
 Assumption \ref{ass:assumption-on-kernel} ({\it 1}) and that $\eps\leq e^{-1}$.

\smallskip

With $\lambda_r(M,\K,Q)=2M+\Lambda(\K,Q) \mathfrak{c}_{2r}$ the use of the Cauchy-Schwartz inequality leads to the following bound:
\begin{eqnarray}
\label{eq8:proof-of-theorem-local}
&&\left\{\mathbb{E}_F^\eps\left|\hat F_{(\hat{\ta},\hat{h})}(x)-F(x)\right|^r\mathrm{1}_{\overline{\mathcal{C}}}\right\}^{1/r}
\lambda_r(M,\K,Q) \| \K \|^2_{\infty}\sqrt{\ln (1/\eps)}
\left[\mathbb{P}_F^\eps(\overline{\mathcal{A}})+\mathbb{P}_F^\eps(\overline{\mathcal{B}})\right]^{1/2r}.
\end{eqnarray}

$1^0.\;$ Let us bound from above $\mathbb{P}_F^\eps(\overline{\mathcal{A}})$. We note that
\begin{eqnarray}
\label{eq9:proof-of-theorem-local}
\mathbb{P}_F^\eps(\overline{\mathcal{A}})
&=&\mathbb{P}_F^\eps\big\{(\ta^{\circ},h^*)\notin \cc P(x)\big\}
= \mathbb{P}_F^\eps\big\{R_{(\ta^{\circ}, h^*)}(x)>0\big\}\nn*[2mm]
&\leq&\sum_{k:\;\eps^2\leq 2^{-k}\leq h^*}\mathbb{P}_F^\eps\left\{\sup_{\nu\in\SS}
\big|\hat F_{(\ta^{\circ}, 2^{-k})(\nu, 2^{-k})}(x) - \hat F_{(\nu, 2^{-k})}(x) \big|>   \TH\big(2^{-k}\big)\right\}.
\end{eqnarray}
For any $k$ satisfying $2^{-k}\leq h^*$ and any $\nu\in\SS$, similarly to (\ref{eq5:proof-of-theorem-local}), we obtain from the first assertion of Lemma  \ref{lem:bounds for bias}
with $h=2^{-k}$, (\ref{eq3:proof-of-theorem-local}) and (\ref{eq4:proof-of-theorem-local}) that
\begin{eqnarray}
\label{eq10:proof-of-theorem-local}
\hskip-0.7cm\left|\hat F_{(\ta^{\circ},2^{-k})(\nu,2^{-k})}(x)-\hat F_{(\nu,2^{-k})}(x)\right|
&\le& 2(h^*_f)^{-1/2}
\| \K \|^{2}_{\infty}\;\eps \sqrt{\ln (1/\eps)}+2\sqrt{2^k}\;\| \K \|^{2}_{2}\eps \sqrt{\ln (1/\eps)}\;\zeta_\eps(x)
\nonumber\\*[2mm]
&\le& 2^{1+k/2}\;\| \K \|^{2}_{\infty}\eps \sqrt{\ln (1/\eps)}\big[1+\zeta_\eps(x)\big].
\end{eqnarray}
Here we have also used that $h^*_f\geq 2^{-k}$. Remembering, that
$$
\TH(\eta)=2\| \K \|^{2}_{\infty}\left[\Lambda(\K,Q)+\sqrt{4r+2}+1\right]\eps\sqrt{\eta^{-1}\ln (1/\eps)},
$$
we obtain from (\ref{eq10:proof-of-theorem-local}) for any $k$ satisfying $2^{-k}\leq h^*$
\begin{eqnarray*}
%\label{eq9:proof-of-theorem-local}
&&\mathbb{P}_F^\eps\left\{\sup_{\nu\in\SS}
\big|\hat F_{(\ta^{\circ}, 2^{-k})(\nu, 2^{-k})}(x) - \hat F_{(\nu, 2^{-k})}(x) \big|>   \TH\big(2^{-k}\big)\right\}
\nonumber
\\*[2mm]
&&\leq \mathbb{P}_F^\eps\bigg\{\displaystyle{\sup_{E\in \mathcal{E}_{\frac{1}{4},\frac{1}{\eps^{2}}}}}\left|\zeta_x\big(E\big)\right|
> \mathbf{c}\big(1/4,\eps^{-2}\big)+ \sqrt{(4r+2)\ln (1/\eps)}\bigg\}\le \eps^{2r+1},
\end{eqnarray*}
in view of the first assertion of Lemma \ref{lem:gauss-on-matrices}. It yields, together with (\ref{eq9:proof-of-theorem-local})
\begin{eqnarray}
\label{eq11:proof-of-theorem-local}
\mathbb{P}_F^\eps(\overline{\mathcal{A}})\le 2 \eps^{2r+1}\log_2 (1/\eps) \le 2 \eps^{2r}.
\end{eqnarray}

$2^0.\;$ An upper bound on the probability of event \( \left\{\hat{h}< h^*\right\} \) is given by
\begin{eqnarray}
\label{eq12:proof-of-theorem-local}
\mathbb{P}_F^\eps(\overline{\mathcal{B}})&=&\mathbb{P}_F^\eps\left(\bigcup_{k:\;\eps^2\leq 2^{-k}\leq h^*}
\left\{\left|\hat F_{(\hat{\ta},h^*)}(x)- \hat F_{(\hat{\ta}, 2^{-k})}(x)\right|
                >\TH\big(2^{-k}\big)\right\}\right)=
\nonumber\\*[2mm]
&\leq&\sum_{k:\;\eps^2\leq 2^{-k}\leq h^*}\mathbb{P}_F^\eps
\left\{\left|\hat F_{(\hat{\ta},h^*)}(x)- \hat F_{(\hat{\ta}, 2^{-k})}(x)\right| >   \TH\big(2^{-k}\big)\right\}.
\end{eqnarray}
We note  that
\begin{eqnarray*}
&&\left|\hat F_{(\hat{\ta},h^*)}(x)-\hat F_{(\hat{\ta},2^{-k})}(x)\right|
\leq \left|S_{(\hat{\ta},h^*)}(x)-S_{(\hat{\ta},2^{-k})}(x)\right|
\\*[2mm]
&&+
 \eps\|K\|_2\sqrt{\det\big( E_{(\hat{\ta},h^*)}\big) }\left|\xi _{(\hat{\ta},h^*)}(x)\right|
+\eps\|K\|_2\sqrt{\det\big( E_{(\hat{\ta},2^{-k})}\big) }\left|\xi _{(\hat{\ta},2^{-k})}(x)\right|.
\end{eqnarray*}
Applying the second assertion of Lemma  \ref{lem:bounds for bias}
with $\nu=\hat{\ta},\;h=h^*,\;\eta=2^{-k}$, (\ref{eq3:proof-of-theorem-local}) and (\ref{eq4:proof-of-theorem-local})
\begin{eqnarray}
\label{eq13:proof-of-theorem-local}
\left|\hat F_{(\hat{\ta},h^*)}(x)-\hat F_{(\hat{\ta},2^{-k})}(x)\right|
&\leq& 2(h^*_f)^{-1/2}\; \| \K \|^{2}_{\infty}\;\eps \sqrt{\ln (1/\eps)}+2\sqrt{2^k}\;\| \K \|^{2}_{2}\eps \sqrt{\ln (1/\eps)}\;\zeta_\eps(x)
\nn *[2mm]
&\leq& 2^{1+k/2}\| \K \|^{2}_{\infty}\eps \sqrt{\ln (1/\eps)}\big[1+\zeta_\eps(x)\big].
\end{eqnarray}
We remark that the right-hand sides of (\ref{eq10:proof-of-theorem-local}) and (\ref{eq13:proof-of-theorem-local}) coincide and, therefore,
repeating the computation led to (\ref{eq11:proof-of-theorem-local}) we get
\begin{eqnarray}
\label{eq14:proof-of-theorem-local}
\mathbb{P}_F^\eps(\overline{\mathcal{B}})\le 2 \eps^{2r}.
\end{eqnarray}
We obtain from (\ref{eq8:proof-of-theorem-local}), (\ref{eq11:proof-of-theorem-local}) and (\ref{eq14:proof-of-theorem-local})
\begin{eqnarray}
\label{eq15:proof-of-theorem-local}
\left\{\mathbb{E}_F^\eps\left|\hat F_{(\hat{\ta},\hat{h})}(x)-F(x)\right|^r\mathrm{1}_{\overline{\mathcal{C}}}\right\}^{1/r}
&\leq& 2^{1/r}\lambda_r(M,\K,Q) \| \K \|^2_{\infty}\eps\sqrt{\ln (1/\eps)}.
%\nonumber\\
%&\leq&C_{r,2}\| \K \|^2_{\infty}(h^*_f)^{-1/2}\eps\sqrt{\ln (1/\eps)},
\end{eqnarray}
%since $h^*_f\leq 1$.
The assertion of the theorem follows now from (\ref{eq7:proof-of-theorem-local}) and (\ref{eq15:proof-of-theorem-local}).
\epr
\subsection{Proof of Theorem \ref{th:pointwise-adaptation}}

We start this section with an auxiliary result used in the proof of the second assertion of the theorem. That result is proved in \cite{KLP2008}, Proposition 7, and for convenience, we formulate it as Lemma \ref{lem:KLP-result} below.

\subsubsection{Auxiliary result}

The result cited below concerns a lower bound for estimators of an arbitrary mapping in the framework of GWN model. Below a version adjusted to the estimation at a given point is provided.

\smallskip

Let \( \cc F \) be a nonempty class of functions and let \( F: \RR^d \to \RR \) be an unknown signal from model \eqref{WGN_model}--\eqref{single-index} satisfying \( F \in \mathcal{F} \subset \mathbb{L}_2(\cc D) \), \( \cc D =[-1, 1]^d \). The aim is to estimate the functional \( F(x) \), \( x \in [-1/2, 1/2]^d \).

\begin{lemma}
\label{lem:KLP-result} (\cite{KLP2008})
Assume that for any \( \eps > 0 \) there exist a positive integer \( \numhyp \), \( \numhyp \to \infty \) as \( \eps \to 0 \),  \( \rho \in (0,1) \),
\( c>0 \) and functions \( F_0, F_1,  \ldots , F_{\numhyp} \in \cc F \) such that:
\begin{eqnarray}
\label{eq:ass1-klp-lemma}
  &&|F_i(x) - F_0(x)| =\lambda_\eps,\qquad\; \forall i=1, \ldots, \numhyp;
\\*[2mm]
\label{eq:ass2-klp-lemma}
&& \langle F_i - F_0 , F_j - F_0 \rangle \le c \eps^2\quad\;\;
                                    \forall i,j=1, \ldots, \numhyp,\;  i \ne j;
\\*[2mm]
\label{eq:ass3-klp-lemma}
&&\| F_i - F_0 \|^2_2 \le \rho \eps^2 \ln (\numhyp),\quad
                                    \forall  i=1, \ldots, \numhyp.
\end{eqnarray}
Then for \( r\geq 1 \)
\begin{equation*}
\label{lower bound GWN}
   \inf_{\tilde F}
        \sup_{F \in \cc F}
            \left(\expect{F}
        \big| \tilde{F} (x) - F (x) \big|^r\right)^{\frac{1}{r}}
   \ge
    \frac{1}{2}\left(  1- \sqrt{\frac{e^c -1}{e^c +3}  }   \right) \lambda_\eps.
   % 2^{-1} \left[ 1-\big(1+4e^{-c}\big)^{-1/2}\right]\lambda_\eps.
\end{equation*}
\end{lemma}

\subsubsection{Proof of Theorem \ref{th:pointwise-adaptation}}
\paragraph{Proof of the first assertion}
Under Assumptions \ref{ass:assumption-on-kernel} and \ref{ass2:assumption-on-kernel} the standard computation of the bias of kernel estimators, for any $f\in\mathbb{H}(\beta,L)$ and any \( z\in\RR \), gives
\begin{equation*}
%\left|
\Delta_{\K,f}(h,z)
%\right|
    \leq
    \frac{ L h^\beta2^{-\beta} \|\K\|_\infty }{ (1+\beta)m_{\beta}! }
    \leq
    \|\K\|_\infty Lh^\beta.
\end{equation*}
The right-hand side of the latter inequality does not depend of \( z \) so
\begin{equation*}
%\left|
\Delta^*_{\K,f}(h,z)
%\right|
    \leq
        \|\K\|_\infty Lh^\beta.
\end{equation*}
Hence,
$h^*_{\K,f}(z)\geq \left(L^{-1}\eps\sqrt{\ln (1/\eps)}\right)^{2/(2\beta+1)}$  for any \( z\in\RR \) and the first assertion of the theorem
follows from Theorem \ref{th:local-oracle-inequality}.

\paragraph{Proof of the second assertion}
The proof is based on the construction of a family $F_0, \ldots,F_{\numhyp}\in \cc F = \mathbb{F}_d(\beta,L) \subset L_2([-1,1]^d)$ satisfying conditions
\eqref{eq:ass1-klp-lemma}--\eqref{eq:ass3-klp-lemma} of Lemma \ref{lem:KLP-result}.

\smallskip

$1^0.\; $ Firstly, we construct \( F_0,\ldots,F_{\numhyp} \)  and verify \eqref{eq:ass1-klp-lemma}.
Let \( g:\RR \to \RR \) be such that
\( \supp (g) \subset (-1/2, 1/2) \),  \( g \in \mathbb{H}(\beta,1) \)%, \( \| g\|_{\infty} \le 1 \),
and $g(0)\neq 0$.
Put  \( h=\left(\mathfrak{a}L^{-1}\eps\sqrt{\ln (1/\eps)}\right)^{2/(2\beta+1)} \),  where the constant $\mathfrak{a}>0$
will be chosen later in order to satisfy \eqref{eq:ass3-klp-lemma}.
For any fixed \( u \in \RR \) define
\begin{equation}\label{link function for hypothesis}
    f_u(v) = Lh^\beta g\big[(v-u) h^{-1}\big] \,, \;\; v\in\RR.
    % f(z) = Lh^\beta g\big[(z-x) h^{-1}\big] \;, \;\; z\in\RR.
\end{equation}
For \( b>0 \) put \( \numhyp = \eps^{-b} \) assuming without loss of generality that \( \numhyp  \) is an  integer.
The value of \( b \) will be determined later in order to satisfy  \eqref{eq:ass2-klp-lemma}.

Let $\left\{\vartheta_i,\; i=1,\ldots, \numhyp \right\}\subset \SSd$ be defined as follows:
\begin{equation*}
    \vartheta_i=\big(\theta^{(1)}_i,\theta^{(2)}_i,0,\ldots,0\big)^{\T},
    \qquad
    \theta^{(1)}_i=\cos(i/\numhyp),
    \quad
    \theta^{(2)}_i=\sin(i/\numhyp).
\end{equation*}
Finally, set
\begin{equation}\label{hypotheses}
    F_0 \equiv 0 \;\; \text{  and }\;\; F_i(t) = f_{\vartheta_i^{\T} x}\big(\vartheta_i^{\T} t\big) \;,\;\; i = 1, \ldots, \numhyp.
\end{equation}
As \( g \in \mathbb{H}(\beta,1) \) so \( f_u \) defined by \eqref{link function for hypothesis} belongs to \( \mathbb{H}(\beta,L) \) for any \( u \in \RR \) and therefore
all \( F_i \) are in \( \cc F= \mathbb{F}_d(\beta,L) \). Moreover, for any $ i=1,\ldots, N_\eps$
\begin{eqnarray}
\label{eq:check-first-condition}
    \big|F_i(x) - F_0(x)\big|
        =
        \big|f_{\vartheta_i^{\T} x}\big(\vartheta_i^{\T} x\big)\big|
        &=&
         |g(0)| L^{\frac{1}{2\beta+1}}
            \left(\mathfrak{a}\eps\sqrt{\ln (1/\eps)}\right)^{\frac{2\beta}{2\beta+1}}\nn
        &=&
          |g(0)| \mathfrak{a}^{\frac{2\beta}{2\beta+1}}\psi_\eps(\beta,L).
\end{eqnarray}
We see that \eqref{eq:ass1-klp-lemma} holds with
\( \lambda_\eps = |g(0)| \mathfrak{a}^{\frac{2\beta}{2\beta+1}}\psi_\eps(\beta,L) \).

\smallskip

$2^0.\; $ Now we  check (\ref{eq:ass2-klp-lemma}).
Set  \( \ort{\theta_{i}}=(-\sin(i/\numhyp)),\cos(i/\numhyp) \). We have
\begin{eqnarray*}
  \langle F_i, F_j \rangle
  &=&
    L^2 h^{2\beta}
        \int_{[-1 , 1 ]^d}
        %\int_{[-\frac{1}{2} , \frac{1}{2} ]^d}
           g\big(h^{-1}\vartheta_i^{\T}(t-x)\big)
           g\big(h^{-1}\vartheta_j^{\T}(t-x)\big)\dd t \\
  &\le&
        3^{d-2}L^2 h^{2\beta+2} \int_{\mathbb{R}^2}\left|g\big(\theta_i^{\T}u\big)
           g\big(\theta_j^{\T}u\big)\right|\dd u
           =
           3^{d-2}L^2 h^{2\beta+2} \big|\ort{\ta_{i}}^{\T} \ta_j\big|^{-1}\| g \|_{1}^2.
\\
&=&
3^{d-2}L^2 h^{2\beta+2} \big|\cos(j/\numhyp)\sin(i/\numhyp)-\cos(i/\numhyp)\sin(j/\numhyp)\big|^{-1}\| g \|_{1}^2
\\*[2mm]
&=&
3^{d-2}L^2 h^{2\beta+2} \big|\sin\big((i-j)/\numhyp\big)\big|^{-1}\| g \|_{1}^2
\\*[2mm]
&=&3^{d-2}L^2 h^{2\beta+2}\big(\sin\big(|i-j|/\numhyp\big)\big)^{-1}\| g \|_{1}^2.
\end{eqnarray*}
Thus, we obtain
\begin{eqnarray}
\label{eq1:proof-the-local-adaptive}
\sup_{i \ne j ; \;i,j = 1, \ldots, \numhyp}  \langle F_i, F_j \rangle
    & \le &
       3^{d-2} L^2 h^{2\beta+2} \big(\sin\big(1/\numhyp\big)\big)^{-1}\| g \|_{1}^2
            \le   3^{d-2} 2L^2 h^{2\beta+2}\numhyp \| g \|_{1}^2
\nonumber\\
&=& 3^{d-2} 2 \| g \|_{1}^2\mathfrak{a}^2\eps^{2}\ln (1/\eps)[\numhyp h].
\end{eqnarray}
Hence, choosing $b<2/(2\beta+1)$ we conclude that (\ref{eq:ass2-klp-lemma}) holds with any given $c>0$
for \( \eps \) small enough.

\smallskip

$3^0.\; $ It remains to verify  \eqref{eq:ass3-klp-lemma}. By changing variables for any $i=1,\ldots, N_\eps$
\begin{eqnarray*}
  \| F_i \|_2^2
   & \le &
   3^{d-1}  \| g\|_2^2 L^2 h^{2\beta+1}
                =
               3^{d-1} \| g \|_{2}^2\mathfrak{a}^2 \eps^{2}\ln (1/\eps)
                =
             3^{d-1}   \| g \|_{2}^2 \mathfrak{a}^2 b^{-1} \eps^{2}\ln \big(N_\eps\big).
\end{eqnarray*}
Here the notation \( \| \cdot \|_2 \) stands for the \( L_2 \) norms on \( [-1, 1]^d \) and \( [-1/2, 1/2] \) correspondingly. Choosing $\mathfrak{a}^2= 3^{- d}b\,  \| g \|_{2}^{-2}$ we see that (\ref{eq:ass3-klp-lemma}) is fulfilled with $\rho=1/3$.

In view of (\ref{eq1:proof-the-local-adaptive}) the constants $c$ from (\ref{eq:ass2-klp-lemma}) and $\mathfrak{a}$ are chosen independently of $L$. Thus, the second assertion of the theorem follows  from Lemma \ref{lem:KLP-result}.
\epr

\subsubsection{Proof of Theorem \ref{th:global-adaptation}}

To prove the theorem we will exploit the ideas developed in \cite{LepMamSpok97}. Moreover, our considerations are, to a great degree, based on the technical result of Lemma \ref{lem:l_p-norm-of-bais} below. Its proof is postponed until Appendix.

\begin{lemma}
\label{lem:l_p-norm-of-bais}
 Grant Assumptions \ref{ass:assumption-on-kernel} and \ref{ass2:assumption-on-kernel}. Then, for any $\mathfrak{p}> 1$,
$0<s\leq \bob$, $\cQ>0$,
$$
\sup_{g\in\mathbb{N}_\mathfrak{p}(s,\cQ)}
\left\|\Delta^*_{\K,g}(h,\cdot)\right\|_\mathfrak{p}\leq 2\cQ h^{s}\|\K\|_\infty
(\tau_\mathfrak{p}+1)\left[2^{s \mathfrak{p}}-1\right]^{-\frac{1}{\mathfrak{p}}},\;\;\forall h>0.
$$
Here $\tau_\mathfrak{p}$ is a depending only of \( \mathfrak{p} \) constant from the $(\mathfrak{p},\mathfrak{p})$-strong maximal inequality.

\end{lemma}

\paragraph{Proof of Theorem \ref{th:global-adaptation}}

It is suffice to  prove the theorem only in the case $r\geq p$. Indeed, remind that the risk $\cc R_{r}^{(\eps)}(\cdot, \cdot)$
is described by the \( L_r \) norm on $[-1/2,1/2]$, therefore
$$
\cc R_{r}^{(\eps)}(\cdot, \cdot)\leq \cc R_{p}^{(\eps)}(\cdot, \cdot), \;\;
r\leq p.
$$
Hence the case $r\leq p$ can be reduced to the case $r=p$.

Yet another observation. In view of embedding of Nikol'skii class \( \mathbb{N}_p(\beta,L) \) in the H\"{o}lder class with parameters \( \beta -1/p \) and \( c L \), \( c>0 \), the assumption \( \beta p > 1 \) provides that
\( f\in\mathbb{F}_M \) and the assumptions of Theorem \ref{th:global-oracle-inequality} are fulfilled. Moreover, in order to obtain the desired the assertion it suffices to bound from above
$
\left\| \sqrt{\frac{\| \K \|^2_{\infty}\eps^{2} \ln (1/\eps)}{h^*_{\K,f}(\cdot)}}\right\|_r.
$

\smallskip

Set $\Gamma_0=\big\{y\in[-1/2,1/2]: \;\; h^*_{\K,f}(y)=1\big\}$ and
$\Gamma_k=\big\{y \in [-1/2,1/2]: \;\; h^*_{\K,f}(y) \in (2^{-k}, 2^{-k+1}] \cap[\eps^2,1]\big\}$ for \( k=1,2,\ldots \; \). Later on, the integration over empty set is supposed to be zero.
 We have
\begin{eqnarray*}
  && \left\|\sqrt{\frac{\|\K \|^2_{\infty}\eps^{2}
 \ln (1/\eps)}{h^*_{\K,f}(\cdot)}}\right\|^r_r
  = \sum_{k\geq 1}\int_{\Gamma_k}\left(\sqrt{\frac{\|\K \|^2_{\infty}\eps^{2}
 \ln (1/\eps)}{h^*_{\K,f}(y)}}\right)^r\dd  y+\int_{\Gamma_0}\left(\sqrt{\frac{\|\K \|^2_{\infty}\eps^{2}
 \ln (1/\eps)}{h^*_{\K,f}(y)}}\right)^r\dd  y.
\end{eqnarray*}
%\begin{equation*}
%\label{eq1:proof-theorem-global-adaptation}
%\left\|\sqrt{\frac{\|\K \|^2_{\infty}\eps^{2}
 %\ln (1/\eps)}{h^*_{\K,f}(\cdot)}}\right\|^r_r
 %=
 %\sum_{k\geq 1}\int_{\Gamma_k}\left(\sqrt{\frac{\|\K \|^2_{\infty}\eps^{2}
 %\ln (1/\eps)}{h^*_{\K,f}(y)}}\right)^r\dd  y+\int_{\Gamma_0}\left(\sqrt{\frac{\|\K \|^2_{\infty}\eps^{2}
 %\ln (1/\eps)}{h^*_{\K,f}(y)}}\right)^r\dd  y.
%\end{equation*}
The definition of $\Gamma_0$ implies
\begin{equation}
\label{eq2:proof-theorem-global-adaptation}
\int_{\Gamma_0}\left|\sqrt{\frac{\|\K \|^2_{\infty}\eps^{2}
 \ln (1/\eps)}{h^*_{\K,f}(y)}}\right|^r\dd  t\leq \left[\|\K \|^2_{\infty}\eps^{2}\ln (1/\eps)\right]^{\frac{r}{2}}.
\end{equation}
Assumption
\ref{ass:assumption-on-kernel} ({\it 2}) implies that $\Delta^{*}_{\K,f}(\cdot,y)$ is continuous on $[\eps^2, 1]$, hence for any $k\geq 1$
%We have  in view of (\ref{eq1:def_oracle bandwidth})
\begin{equation}
\label{eq3:proof-theorem-global-adaptation}
\Delta^*_{\K,f}\big( h^*_{\K,f}(y),y\big)=\left[\frac{\| \K \|^2_{\infty}\eps^2 \ln (1/\eps)}{ h^*_{\K,f}(y)}\right]^{\frac{1}{2}},
\;\;\forall y\in\Gamma_k.
\end{equation}
Let $0\leq q_k\leq r$ be a sequence whose choice will be done later. We obtain from (\ref{eq3:proof-theorem-global-adaptation})
\begin{eqnarray}
\label{eq4:proof-theorem-global-adaptation}
\sum_{k\geq 1}\int_{\Gamma_k}\left(\sqrt{\frac{\| \K \|^2_{\infty}\eps^{2}
 \ln (1/\eps)}{h^*_{\K,f}(y)}}\right)^r\dd  y
&\le& \sum_{k\geq 1}\left(\frac{\| \K \|^2_{\infty}\eps^{2}  \ln (1/\eps)}{2^{-k}}\right)^{\frac{r-q_k}{2}}
        \int_{\Gamma_k}\left(\Delta^*_{\K,f}\big(2^{1-k} ,y\big)\right)^{q_k}\dd  y \nn
&\le& \sum_{k\geq 1}\left(\frac{\| \K \|^2_{\infty}\eps^{2}  \ln (1/\eps)}{2^{-k}}\right)^{\frac{r-q_k}{2}}
        \int\left( \Delta^*_{\K,f}\big(2^{1-k} ,y\big)\right)^{q_k}\dd  y=:\Xi.
\end{eqnarray}
To get the first inequality we have used that $\Delta^*_{\K,f}\big(\cdot,y\big)$ in monotonically increasing function.

\par The computation of the quantity on the right-hand side of (\ref{eq4:proof-theorem-global-adaptation}), including the choice of
$(q_k,\;k\geq 1)$,
 will be done differently in dependence on
$\beta, p$ and $r$. Later on $c_1,c_2,\ldots,$ denote constants independent on $\eps$, $L$ and $\K$.

\smallskip

$1^0.\;$ {\it Case $(2\beta+1)p > r$.} Put
$$
h^*=\left[L^{-2}\eps^{2}
 \ln (1/\eps)\right]^{\frac{1}{2\beta+1}}
$$
and choose $q_k=p$ if $2^{-k}\leq h^*$ and $q_k=0$ if $2^{-k}> h^*$.

Applying  Lemma \ref{lem:l_p-norm-of-bais} with $\mathfrak{p}=p$, $s=\beta$ and $\cQ=L$ we get
\begin{eqnarray}
\label{eq55:proof-theorem-global-adaptation}
%&&\hskip-0.7cm
    \Xi
    &\le&  c_1\big(L\|\K\|_\infty\big)^{p}\sum_{k:\;2^{-k}\leq h^*}\left(\frac{\| \K \|^2_{\infty}\eps^{2}
 \ln (1/\eps)}{2^{-k}}\right)^{\frac{r-p}{2}}2^{-k\beta p}+c_2\left(\frac{\| \K \|^2_{\infty}\eps^{2}
 \ln (1/\eps)}{h^*}\right)^{\frac{r}{2}}    \nonumber\\
%&&\hskip-0.7cm
        &\le& c_3 \|\K\|^{r}_\infty\left[L^{p}\left(\eps^{2}
 \ln (1/\eps)\right)^{\frac{r- p}{2}}\sum_{k:\;2^{-k}\leq h^*}2^{-k\big[\beta p-\frac{r-p}{2}\big]}
+\left(\frac{\eps^{2}
 \ln (1/\eps)}{h^*}\right)^{\frac{r}{2}}\right].
\end{eqnarray}
Because in the considered case $\beta p-\frac{r-p}{2}>0$, we obtain
\begin{eqnarray*}
%\label{eq5:proof-theorem-global-adaptation}
\Xi\leq c_4  \|\K\|^{r}_\infty\left[L^{ p}\left(\eps^{2}
 \ln (1/\eps)\right)^{\frac{r- p}{2}}(h^*)^{\beta p-\frac{r- p}{2}}
+\left(\frac{\eps^{2}
 \ln (1/\eps)}{h^*}\right)^{\frac{r}{2}}\right].
\end{eqnarray*}
It remains to note that $h^*$ is chosen by balancing two terms on the right-hand side of the latter inequality. It yields
\begin{eqnarray}
\label{eq5:proof-theorem-global-adaptation}
\Xi\leq 2c_4 \left[\| \K \|_{\infty}L^{\frac{1}{2\beta+1}}\left(\eps\sqrt{\ln(1/\eps)}\right)^{\frac{2\beta}{2\beta+1}}\right]^r.
\end{eqnarray}
The argument in the case $(2\beta+1)p> r$ is completed with the use of Theorem~\ref{th:global-oracle-inequality},
(\ref{eq2:proof-theorem-global-adaptation}) and (\ref{eq5:proof-theorem-global-adaptation}).

\medskip

$2^0.\;$ {\it Case $(2\beta+1)p=r$.} Put $h^*=1$ and choose $q_k=p$ for all $k\geq 1$. Repeating the computations
led to  (\ref{eq55:proof-theorem-global-adaptation}) we get
\begin{eqnarray}
\label{eq6:proof-theorem-global-adaptation}
\Xi\leq c_5\ln(1/\eps)\left[\|\K\|_\infty L^{ p/r}\left(\eps^{2}
 \ln (1/\eps)\right)^{\frac{r-p}{2r}}\right]^r.
\end{eqnarray}
Here we have used that $\beta p-\frac{r-p}{2}=0$ and that the summation in (\ref{eq4:proof-theorem-global-adaptation})
runs over $k$ such that $2^{-k}\geq \eps^{2}$, since otherwise $\Gamma_k=\emptyset$. It remains to note that the equality $(2\beta+1)p=r$
is equivalent to $p/r=1/(2\beta+1)$ and $(r-p)/2r=\beta/(2\beta+1)$. The assertion of the theorem in the case $(2\beta+1)p=r$
follows now from Theorem \ref{th:global-oracle-inequality},
(\ref{eq2:proof-theorem-global-adaptation}) and (\ref{eq6:proof-theorem-global-adaptation}).

\medskip

$3^0.\;$ {\it Case $(2\beta+1)p<r$.} Choose $q_k=r$ if $2^{-k}\leq h^*$ and $q_k=p$ if $2^{-k}> h^*$,
where the choice of $h^*$ will be done later.

The following embedding holds, see  \cite{BesovIlNik}:
$\mathbb{N}_p(\beta,L)\subseteq\mathbb{N}_r\big(\beta-1/p+1/r,c_6L\big)$. Thus, applying Lemma \ref{lem:l_p-norm-of-bais}
with $\mathfrak{p}=r$, $s=\beta-1/p+1/r$ and $\cQ=c_6L$ we get
\begin{eqnarray}
\label{eq7:proof-theorem-global-adaptation}
\Xi_1&:=&\sum_{k:\: 2^{-k}\leq h^*}\left(\frac{\| \K \|^2_{\infty}\eps^{2}
 \ln (1/\eps)}{2^{-k}}\right)^{\frac{r-q_k}{2}}\int\left(
\overline{\Delta}_{\K,f}\big(2^{1-k} ,y\big)\right)^{q_k}\dd  y
\nonumber\\
&=&\sum_{k:\: 2^{-k}\leq h^*}\int\left(
\overline{\Delta}_{\K,f}\big(2^{1-k} ,y\big)\right)^{r}\dd  y
\leq c_7 \big(\|\K\|_\infty L\big)^{r}(h^*)^{\beta r -(r/p)+1}.
\end{eqnarray}
Applying Lemma \ref{lem:l_p-norm-of-bais}
with $\mathfrak{p}=r$, $s=\beta$ and $\cQ=L$ we get
\begin{eqnarray}
\label{eq8:proof-theorem-global-adaptation}
\Xi_2&:=&\sum_{k:\: 2^{-k}> h^*}\left(\frac{\| \K \|^2_{\infty}\eps^{2}
 \ln (1/\eps)}{2^{-k}}\right)^{\frac{r-q_k}{2}}\int\left(
\overline{\Delta}_{\K,f}\big(2^{1-k} ,y\big)\right)^{q_k}\dd  y
\nonumber\\
&=&c_8L^p\big(\|\K\|_\infty\big)^{r}\left(\eps^{2}
 \ln (1/\eps)\right)^{\frac{r-p}{2}}\sum_{k:\: 2^{-k}> h^*}2^{-k\big[\beta p-\frac{r-p}{2}\big]}
\nonumber\\
&\leq&c_9L^p\big(\|\K\|_\infty\big)^{r}\left(\eps^{2}
 \ln (1/\eps)\right)^{\frac{r-p}{2}}(h^*)^{\beta p-\frac{r-p}{2}}.
\end{eqnarray}
Here we have used that $\beta p-\frac{r-p}{2}<0.$
In view of (\ref{eq7:proof-theorem-global-adaptation}) and (\ref{eq8:proof-theorem-global-adaptation}) we choose $h^*$ from the equality:
$$
L^{r}(h^*)^{\beta r -(r/p)+1}=L^p\left(\eps^{2}
 \ln (1/\eps)\right)^{\frac{r-p}{2}}(h^*)^{\beta p-\frac{r-p}{2}}.
$$
It yields $h^*=\left(L^{-1}\eps\sqrt{\ln (1/\eps)}\right)^{\frac{1}{\beta-1/p+1/2}}$ and we obtain finally that
\begin{eqnarray}
\label{eq9:proof-theorem-global-adaptation}
\Xi\leq c_{10}\big(\|\K\|_\infty\big)^{r}L^{\frac{r(1/2-1/r)}{\beta-1/p+1/2}}
\left(\eps\sqrt{\ln (1/\eps)}\right)^{\frac{r(\beta-1/p+1/r)}{\beta-1/p+1/2}}.
\end{eqnarray}
The assertion of the theorem in the case $(2\beta+1)p<r$
follows now from Theorem \ref{th:global-oracle-inequality},
(\ref{eq2:proof-theorem-global-adaptation}) and (\ref{eq9:proof-theorem-global-adaptation}). \epr
\section{Appendix}\label{sec:appendix}
\subsection{Proof of Lemma \ref{lem:bounds for bias}}
\paragraph{Proof of the first assertion}
The symmetry of the kernel $\K$ (Assumption~\ref{ass:assumption-on-kernel}~({\it 1})) implies
$$
S_{(-\ta^{\circ}, h)(\nu, h)}(\cdot)\equiv S_{(\ta^{\circ}, h)(-\nu, h)}(\cdot),\quad S_{(-\nu, h)}(\cdot)\equiv S_{(\nu, h)}(\cdot).
$$
Therefore it suffices to prove the first assertion of the lemma under the condition $\nu^\T\ta^{\circ}\geq 0$. In this case
$E_{(\ta^{\circ}, h)(\nu, h)}=\overline{E}_{(\ta^{\circ}, h)(\nu, h)}$ and we note that
%The following relation  will be useful in the sequel: for any $\ta,\nu\in\SS$ and any $h\in[\eps^2,1]$
\begin{equation}
\label{eq1:matrix for pseudo-estimator}
\overline{E}_{(\ta^{\circ}, h)(\nu, h)}
=
\left[E^{-1}_{(\ta^{\circ}, h)} + E^{-1}_{(\nu, h)}\right]^{-1}.
\end{equation}
For any $\theta=(\theta_1,\theta_2)\in\SS$ let $\ort{\ta}=(-\theta_2,\theta_1)$.
Using (\ref{eq1:matrix for pseudo-estimator}) we obtain
\begin{eqnarray*}
%\label{pseudo_smooth_afterCHV}
S_{(\ta^{\circ}, h)(\nu, h)}(x)
    =
\int K(u)
    f\big(h [\ta^{\circ} + \nu]^{\T} \ta^{\circ} u_1 + [\ort{\ta^{\circ}} + \ort{\nu}]^{\T} \ta^{\circ} u_2
        + x^{\T} \ta^{\circ}\big) \dd u
\\
=
\int\int \K(u_1)\K(u_2)f\big( h [1 + \nu^{\T} \ta^{\circ}] u_1 + \ort{\nu}^{\T} \ta^{\circ} u_2+ x^{\T} \ta^{\circ}\big) \dd u_1\dd u_2.
\end{eqnarray*}
We also have
$$
S_{(\nu, h)}(x)=
\int\int \K(u_1)\K(u_2)f\big( h \nu^{\T} \ta^{\circ} u_1 + \ort{\nu}^{\T} \ta^{\circ} u_2+ x^{\T} \ta^{\circ}\big) \dd u_1\dd u_2.
$$
Put
$
S_\nu^*(x)=\int \K(u_2)f\big(\ort{\nu}^{\T} \ta^{\circ} u_2+ x^{\T} \ta^{\circ}\big) \dd u_2
$
and consider two cases.

\medskip

$1^{0}.\;\; \ort{\nu}^{\T} \ta^{\circ}=0$. In this case $S_\nu^*(x)=f(x^{\T} \ta^{\circ})$ and
\begin{gather*}
S_{(\nu, h)}(x)=\int \K(u_1)f\big( h u_1 + x^{\T} \ta^{\circ}\big) \dd u_1=h^{-1}\int \K\big([t-x^{\T} \ta^{\circ}]/h\big)f(t) \dd t,
\\
S_{(\ta^{\circ}, h)(\nu, h)}(x)=\int \K(u_1)f\big( 2h u_1 +x^{\T} \ta^{\circ}\big) \dd u_1=(2h)^{-1}\int \K\big([t-x^{\T} \ta^{\circ}]/2h\big)f(t) \dd t.
\end{gather*}
Here we have used that  $\ort{\nu}^{\T} \ta^{\circ}=0$ together with $\nu^{\T} \ta^{\circ}\geq 0$ implies $\nu=\ta^{\circ}$.
Thus, we obtain
\begin{eqnarray}
\label{eq2:proof-lemma-bais}
\left|S_{(\ta^{\circ}, h)(\nu, h)}(x)-S_{(\nu, h)}(x)\right|
&\le& \left|S_{(\ta^{\circ}, h)(\nu, h)}(x)-S_\nu^*(x)\right|+\left|S_{(\nu, h)}(x)-S_\nu^*(x)\right|
\nonumber\\*[2mm]
&\le&\Delta_{\K,f}\big(h,x^{\T} \ta^{\circ}\big)+\Delta_{\K,f}\big(2h,x^{\T} \ta^{\circ}\big)\leq 2\Delta^*_{\K,f}\big(2h,x^{\T} \ta^{\circ}\big).
\end{eqnarray}

\smallskip

$2^{0}.\;\; \ort{\nu}^{\T} \ta^{\circ}\neq0$. In this case we have
\begin{eqnarray*}
&&S_\nu^*(x)
=\int\int\frac{1}{h(1+\nu^{\T} \ta^{\circ})}\K\left(\frac{v_1}{h(1+\nu^{\T} \ta^{\circ})}\right)
\frac{1}{|\ort{\nu}^{\T} \ta^{\circ}|}\K\left(\frac{v_2-x^{\T} \ta^{\circ}}{|\ort{\nu}^{\T} \ta^{\circ}|}\right)f(v_2) \dd v_1\dd v_2,
\\*[2mm]
%S_{(\nu, h)}(x)=\int \K(u_1)f\big(\pm h u_1 + x^{\T} \ta^{\circ}\big) \dd u_1=h^{-1}\int \K\big([t-x^{\T} \ta^{\circ}]/h\big)f(t) \dd t;
%\\
&&S_{(\vartheta^*, h)(\nu, h)}(x)
%\\&&
=\int \int \frac{1}{h(1+\nu^{\T} \ta^{\circ})}\K\left(\frac{v_1}{h(1+\nu^{\T} \ta^{\circ})}\right)
\frac{1}{|\ort{\nu}^{\T} \ta^{\circ}|}\K\left(\frac{v_2-x^{\T} \ta^{\circ}}{|\ort{\nu}^{\T} \ta^{\circ}|}\right)f(v_1+v_2) \dd v_1\dd v_2.
\end{eqnarray*}
Here we have used once again the symmetry of kernel $\K$. Thus, taking into account that $|\nu^{\T} \ta^{\circ}|\leq 1$, we get
\begin{eqnarray*}
&&\left|S_{(\ta^{\circ}, h)(\nu, h)}(x)-S_\nu^*(x)\right|
\\
&&\leq\int \frac{1}{|\ort{\nu}^{\T} \ta^{\circ}|}\left|\K\left(\frac{v_2-x^{\T} \ta^{\circ}}{|\ort{\nu}^{\T} \ta^{\circ}|}\right)\right|
 \sup_{\delta\leq 2h}\left|\int\frac{1}{\delta}\K\left(\frac{v_1}{\delta}\right)
\big[f(v_1+v_2)-f(v_2)\big] \dd v_1\right|\dd v_2
\\*[2mm]
&&\leq\|\K\|_\infty\sup_{a>0}\left[\frac{1}{a}\int_{x^{\T} \ta^{\circ}-a/2}^{x^{\T} \ta^{\circ}+a/2} \sup_{\delta\leq 2h}\left|\int\frac{1}{\delta}\K\left(\frac{v_1}{\delta}\right)
\big[f(v_1+v_2)-f(v_2)\big] \dd v_1\right|\dd v_2\right].
\end{eqnarray*}
Here we have used that $\text{supp}(\K)\subseteq [-1/2,1/2]$ (Assumption \ref{ass:assumption-on-kernel} ({\it 1})). Hence,
\begin{eqnarray}
\label{eq3:proof-lemma-bais}
&&\left|S_{(\ta^{\circ}, h)(\nu, h)}(x)-S_\nu^*(x)\right|
\le \|\K\|_\infty\overline{\Delta}_{\K,f}\big(2h,x^{\T} \ta^{\circ}\big)\leq
\|\K\|_\infty\Delta^*_{\K,f}\big(2h,x^{\T} \ta^{\circ}\big).
\end{eqnarray}
If $\nu^{\T} \ta^{\circ}\neq 0$ we obtain by the same computations
$$
\left|S_{(\nu, h)}(x)-S_\nu^*(x)\right|\leq
\|\K\|_\infty\Delta^*_{\K,f}\big(h,x^{\T} \ta^{\circ}\big).
$$
Noting that $S_{(\nu, h)}(\cdot)\equiv S_\nu^*(\cdot)$ if $\nu^{\T} \ta^{\circ}= 0$ we get
\begin{eqnarray}
\label{eq33:proof-lemma-bais}
\left|S_{(\nu, h)}(x)-S_\nu^*(x)\right|\leq \|\K\|_\infty\Delta^*_{\K,f}\big(h,x^{\T} \ta^{\circ}\big),
\end{eqnarray}
that yields together with (\ref{eq3:proof-lemma-bais})
\begin{eqnarray}
\label{eq4:proof-lemma-bais}
\left|S_{(\ta^{\circ}, h)(\nu, h)}(x)-S_{(\nu, h)}(x)\right|\leq
2\|\K\|_\infty\Delta^*_{\K,f}\big(2h,x^{\T} \ta^{\circ}\big).
\end{eqnarray}
Finally, taking into account that in view of Assumption \ref{ass:assumption-on-kernel} ({\it 1}) $\|\K\|_\infty\geq 1$, we obtain
from (\ref{eq2:proof-lemma-bais}) and (\ref{eq4:proof-lemma-bais}) that
$$
\left|S_{(\ta^{\circ}, h)(\nu, h)}(x)-S_{(\nu, h)}(x)\right|\leq
2\|\K\|_\infty\Delta^*_{\K,f}\big(2h,x^{\T} \ta^{\circ}\big)\leq 2\|\K\|_\infty\Delta^*_{\K,f}\big(h^*_f,x^{\T} \ta^{\circ}\big),
$$
since we consider $h$ such that $2h\leq h^*_{f}$. The definition of $h^*_{f}$ implies
$$
\Delta^*_{\K,f}\big(h^*_f,x^{\T} \ta^{\circ}\big)\leq (h^*_f)^{-1/2} \| \K \|_{\infty}\eps \sqrt{\ln (1/\eps)}
$$
and the first assertion of the lemma follows.
%since $\|\K\|_\infty\geq \|\K\|_2$ in view of $\text{supp}(\K)\subseteq [-1/2,1/2]$.

\paragraph{Proof of the second and third  assertions} In view of (\ref{eq33:proof-lemma-bais}) for  $\forall\eta\leq h\leq  h^*_{f}$
\begin{eqnarray*}
\left|S_{(\nu, \eta)}(x)-S_{(\nu, h)}(x)\right|&\leq&\left|S_{(\nu, \eta)}(x)-S_\nu^*(x)\right|+\left|S_{(\nu, h)}(x)-S_\nu^*(x)\right|
\\
&\le&\|\K\|_\infty\left[\Delta^*_{\K,f}\big(\eta,x^{\T} \ta^{\circ}\big)+\Delta^*_{\K,f}\big(h,x^{\T} \ta^{\circ}\big)\right]
\le
2\|\K\|_\infty\Delta^*_{\K,f}\big(h,x^{\T} \ta^{\circ}\big)
\\
&\leq& 2\|\K\|_\infty\Delta^*_{\K,f}\big(h^*_f,x^{\T} \ta^{\circ}\big)\leq 2(h^{*}_f)^{-1/2} \| \K \|^2_{\infty}\eps \sqrt{\ln (1/\eps)},
\end{eqnarray*}
in view of the definition of $ h^*_{f}$. The second assertion is proved.

\smallskip

%We remark that $S_{(\ta^{\circ}, h)} (x)=S_{(-\ta^{\circ}, h)} (x)$ since the kernel $\K$ is symmetric. Moreover,
We have for any $ h\leq  h^*_{f}$
\begin{eqnarray*}
\big|S_{(\ta^{\circ}, h)} (x) - F(x)\big|
    &=&
   \left| \frac{1}{h}\int \K \left( \frac{u}{h}\right) \big[ f(u +x^{\T} \ta^{\circ}) - f(x^{\T} \ta^{\circ})\big] \dd u\right| \\%*[2mm]
   &\le&\Delta_{\K,f}\big(h,x^{\T} \ta^{\circ}\big)
        \le  \Delta^*_{\K,f}\big(h,x^{\T} \ta^{\circ}\big)\leq  \Delta^*_{\K,f}\big(h^*_f,x^{\T} \ta^{\circ}\big)
        \\%*[2mm]
        &=&(h^*_f)^{-1/2} \| \K \|_{\infty}\eps \sqrt{\ln (1/\eps)},
\end{eqnarray*}
in view of the definition of $ h^*_{f}$. The third assertion is proved. \epr

\subsection{Proof of Lemma \ref{lem:gauss-on-matrices}}
Since $\zeta_{y}(\cdot)$ is a zero mean Gaussian random function
%in view of the obvious relation $\sup_E|\zeta_{y}(E)|=[\sup_E \zeta_{y}(E)]\vee [\sup_E \{-\zeta_{y}(E)\}]$
we have
\begin{equation}
\label{eq1:proof-lemma-gauss-on-matrices}
\PP\left\{\sup_{E\in\mathcal{E}_{a,A}}\left|\zeta_{y}(E)\right|\geq u \right\}
\leq
2
\PP\left\{\sup_{E\in\mathcal{E}_{a,A}}\zeta_{y}(E)\geq u \right\},\quad\forall u>0.
\end{equation}
%Let \( \EE\Big[\sup_{E\in\mathcal{E}_{a,A}}\zeta_{y}(E)\Big] \le \mathbf{c}(a,A)\), where the value of \( \mathbf{c}(a,A) \) will be determined later.
By Lemma 12.2 in \cite{Lif} the median \( m \) of the random variable \( \sup_{E\in\mathcal{E}_{a,A}}\zeta_{y}(E) \) is dominated by the expectation, that is \( m \le \EE \sup_{E\in\mathcal{E}_{a,A}}\zeta_{y}(E)  \). That along with the Borell, Tsirelson, Sudakov concentration inequality, see Theorem 12.2 in \cite{Lif} provides
\begin{equation}\label{eq02:proof-lemma-gauss-on-matrices}
    \PP\left\{\sup_{E\in\mathcal{E}_{a,A}}\zeta_{y}(E)\ge \EE \sup_{E\in\mathcal{E}_{a,A}}\zeta_{y}(E)  + z \right\}
     \le  \PP\left\{\sup_{E\in\mathcal{E}_{a,A}}\zeta_{y}(E)\ge m + z \right\}
     \le \PP\left\{ \varsigma \ge z\right \}
\end{equation}
%\begin{eqnarray*}
  %\PP\left\{\sup_{E\in\mathcal{E}_{a,A}}\zeta_{y}(E)\ge \mathbf{c}(a,A) + z \right\}
 % &\le&
    %    \PP\left\{\sup_{E\in\mathcal{E}_{a,A}}\zeta_{y}(E)\ge \EE\Big[\sup_{E\in\mathcal{E}_{a,A}}\zeta_{y}(E)\Big] + z \right\}  \\
  %\le  \PP\left\{\sup_{E\in\mathcal{E}_{a,A}}\zeta_{y}(E)\ge m + z \right\}
  %&\le&
     %   \PP\left\{ \varsigma \ge z\right \},
%\end{eqnarray*}
since \( \sup_{E\in\mathcal{E}_{a,A}} \Var \left[  \zeta_{y}(E)  \right] =1\). Here \( \varsigma\;\sim\;\norm{0}{1} \).
%Then using the standard bound for the Gaussian tail we have
%\begin{equation}
%\label{eq11:proof-lemma-gauss-on-matrices}
%\PP\left\{\sup_{E\in\mathcal{E}_{a,A}}\left|\zeta_{y}(E)\right|\ge \mathbf{c}(a,A) + z \right\}\leq
%2 \PP\left\{ \varsigma \ge z\right \} \le \exp\{ -z^2/2 \} ,\quad \forall z>0.
%\end{equation}
Thus, to complete the proof of the first assertion of the lemma it suffices to bound \( \EE \sup_{E\in\mathcal{E}_{a,A}}\zeta_{y}(E) \). This will be done by the application of Dudley's theorem, see Theorem 14.1 in~\cite{Lif}.
Denote by \( \varrho \) the semi-metric generated by \( \zeta_{y}(\cdot) \) on \( \mathcal{E}_{a,A} \):
\begin{equation*}
\varrho(E,E^\prime)=\sqrt{\EE\left|\zeta_{y}(E)-\zeta_{y}(E^\prime)\right|^{2}}, \quad E,E^\prime\in\mathcal{E}_{a,A}.
\end{equation*}
Without loss of generality one can assume that \( \left| \det(E) \right| \ge \left|\det(E^{\prime})\right| \), then we have
\begin{eqnarray*}
&&\varrho^2(E,E^\prime)
=
2\left[1-\| \cc L\|_2^{-2}\sqrt{\left|\det(E)\right|\left|\det(E^{\prime})\right|}
\int  \cc L(Ev)  \cc L(E^{\prime}v)\dd v\right].
\\
&=&2\left[1-
\| \cc L\|_2^{-2}\sqrt{\frac{\left|\det(E^{\prime})\right|}{\left|\det(E)\right|}}
\int  \cc L(z)  \cc L(E^{\prime}E^{-1}z)\dd z\right]
\\
&=&2\left[1-\| \cc L\|_2^{-2}\sqrt{\frac{\left|\det(E^{\prime})\right|}{\left|\det(E)\right|}}
\int_{[-\frac{1}{2} , \frac{1}{2}]^{2}}  \cc L(z)  \cc L\big(E^{\prime}E^{-1}z\big)\dd z\right]
\\
&=&2\left[1- \frac{\sqrt{|\det (E^{\prime})|}}{\sqrt{|\det (E)|}} \right]
+
\frac{2}{\| \cc L \|^2_2} \frac{\sqrt{|\det (E^{\prime})|}}{\sqrt{|\det (E)|}}
\left[
\int_{[-\frac{1}{2} , \frac{1}{2}]^{2}} \cc{L}(z) \big[\cc{L}(z) - \cc{L}(E^{\prime}E^{-1}z)\big]\dd z
 \right].
\end{eqnarray*}
One bounds the first summand with the use of  the assumption \( |\det (E)| \ge a \):
\begin{equation*}
    2\left[1- \frac{\sqrt{|\det (E^{\prime})|}}{\sqrt{|\det (E)|}} \right]
    \le
    \frac{2}{\sqrt{a} }\, \left|\sqrt{|\det (E)|} - \sqrt{|\det (E^{\prime})|}\right|
    \le
     \frac{2}{\sqrt{a} }\, \left|\det (E) - \det (E^{\prime})\right|^{1/2}.
\end{equation*}
As for the second term, putting
$$
\mathfrak{d}^2(E,E^\prime)
=
\int_{[-\frac{1}{2} , \frac{1}{2}]^{2}} \left| \cc L \big(E^{\prime}E^{-1}z\big)- \cc L(z)\right|^{2}\dd z,
$$
by the Cauchy-Schwartz inequality we get
$$
  % \frac{2}{\| \cc L \|^2_2} \frac{\sqrt{|\det (E^{\prime})|}}{\sqrt{|\det (E)|}}
%\left[
\int_{[-\frac{1}{2} , \frac{1}{2}]^{2}} \cc{L}(z) \left[ \cc L (z) - \cc L (E^{\prime}E^{-1}z) \right] \dd z
 %\right]
  \le
 % \frac{2}{\| \cc L \|_2} \frac{\sqrt{|\det (E^{\prime})|}}{\sqrt{|\det (E)|}}
 \| \cc L \|_2 \mathfrak{d}(E,E^\prime).
$$
As $\|\cc L\|_2\geq 1$, we have
$$
\varrho^2(E,E^\prime)\leq 2a^{-1/2}
            \left| \det (E)- \det (E^{\prime}) \right|^{1/2} + 2\mathfrak{d}(E,E^\prime).
$$
First, we note that
$$
\left|\det (E) -\det (E^{\prime})\right|\leq 4A\big|E-E^\prime\big|_\infty.
$$
Second, because \( \cc L \) satisfies the Lipschitz condition with a constant \( \Upsilon \), we have
$$
\mathfrak{d}(E,E^\prime)
    \leq
    \Upsilon  \sup_{z\in[-\frac{1}{2} , \frac{1}{2}]^{2}}\left|(E^\prime-E)E^{-1}z\right|_2
    \leq
     2 \sqrt{2} \Upsilon Aa^{-1} \big|E-E^\prime\big|_\infty.
$$
Since  we assumed  $a\leq A$, the following bound holds:
\begin{equation}
\label{eq2:proof-lemma-gauss-on-matrices}
\varrho^2(E,E^\prime)
\leq
     4 \big(\sqrt{2}\Upsilon+1\big)Aa^{-1}
     \left(\big| E-E^\prime\big|^{1/2}_\infty \bigvee \big|E-E^\prime\big|_\infty \right).
\end{equation}
Consider the cube $[0,A]^{4}$ endowed with the vector supremum norm \( |z|_\infty=\max_{i=1,\ldots, 4}|z_i| \). Let
\( \mathfrak{E}_{[0,A]^{4},|\cdot|_\infty}(\cdot) \) denote the metric entropy of \( [0,A]^{4} \) measured in \( |\cdot|_\infty \). Then
\begin{equation*}
    \mathfrak{E}_{[0,A]^{4},|\cdot|_\infty}(\epsilon)\leq 4\ln(A)+\left[  4\ln{(1/(2\epsilon))}  \right]_{+},\quad \forall\epsilon\in (0,1].
\end{equation*}
Denoting by \( \mathfrak{E}_{\mathcal{E}_{a,A},\varrho}(\cdot) \) the metric entropy of
\( \mathcal{E}_{a,A} \) measured in \( \varrho \), we get in view of (\ref{eq2:proof-lemma-gauss-on-matrices})
\begin{equation*}
    \mathfrak{E}_{\mathcal{E}_{a,A},\varrho}(\delta)
\leq
\mathfrak{E}_{[0,A]^{4},|\cdot|_\infty}
        \left(
            \frac{\delta^{4}}{16(1 + \sqrt{2} \Upsilon)^2A^2a^{-2}}
        \right),\quad\forall\delta\in (0,1],
\end{equation*}
and, therefore,
\begin{equation*}
    \mathfrak{E}_{\mathcal{E}_{a,A},\varrho}(\delta)
\leq
4 \left[ \ln\big( A\vee\{A/a\}^2 \big) + \ln 8 + 2\ln{(1+\sqrt{2}\Upsilon)} +4\ln{(1/\delta)}\right].
\end{equation*}
Since \(  \sup_{E\in\mathcal{E}_{a,A}} \Var \left[ \zeta_{y}(E)  \right] =1\) the use of Dudley's integral bound, see Theorem 14.1 in~\cite{Lif}, leads to
\begin{eqnarray}
\label{eq100:proof-lemma-gauss-on-matrices}
 \EE\Big[\sup_{E\in\mathcal{E}_{a,A}}\zeta_{y}(E)\Big]
&\le&
    4\sqrt{2}\int_{0}^{1/2}\sqrt{\mathfrak{E}_{\mathcal{E}_{a,A},\varrho}(\delta)}\;\dd\delta
\nonumber\\
&\le&
4\sqrt{2}  \left[ \ln(A\vee\{A/a\}^2) + 2 \ln{(1+ \sqrt{2} \Upsilon)} \right]^{1/2} + 29
       \; =:  \, \mathbf{c}(a,A).
\end{eqnarray}
Here we have used that \( \int_{0}^{1/2}\sqrt{\ln{(1/\delta)}}\;\dd\delta\leq 2^{-1}\sqrt{\pi} \).
The first assertion of the lemma follows now from \eqref{eq1:proof-lemma-gauss-on-matrices}, \eqref{eq02:proof-lemma-gauss-on-matrices}, \eqref{eq100:proof-lemma-gauss-on-matrices} and the standard bound for the Gaussian tail.
 %and the Borel-Tsirelson-Sudakov concentration inequality, see e.g. Proposition 2.2 in \cite{Talagrand}.

\smallskip

To justify the second assertion we first note that for any \( q\ge1 \)
$$
\mathbb{E}\Big[\sup_{E\in\mathcal{E}_{a,A}}\big|\zeta_{y}(E)\big|\Big]^{q}=
q\int_{0}^{\infty}u^{q-1}\mathbb{P}\left(\sup_{E\in\mathcal{E}_{a,A}}\big|\zeta_{y}(E)\big|\geq u\right)\dd u.
$$
Hence, applying the first assertion of the lemma we have
\begin{eqnarray*}
% \nonumber to remove numbering (before each equation)
\mathbb{E}\Big[\sup_{E\in\mathcal{E}_{a,A}}\big|\zeta_{y}(E)\big|\Big]^{q}
      &\le&  \big[\mathbf{c}(a,A)\big]^{q} + q \int_0^{\infty}\PP\left\{ |\varsigma| \ge z\right \}(\mathbf{c}(a,A) + z)^{q-1} \dd z\\
   &=& \mathbb{E}\big(\mathbf{c}(a,A)+|\varsigma|\big)^{q},
\end{eqnarray*}
%$$
%\mathbb{E}\Big[\sup_{E\in\mathcal{E}_{a,A}}\big|\zeta_{y}(E)\big|\Big]^{q}\leq \big[\mathbf{c}(a,A)\big]^{q}+\sqrt{2\pi}\;\mathbb{E}\big(\mathbf{c}(a,A)+|\varsigma|\big)^{q},
%$$
where \( \varsigma\;\sim\;\norm{0}{1}  \). Thus, we finally have
\begin{equation*}
    \left(\mathbb{E}\Big[\sup_{E\in\mathcal{E}_{a,A}}\big|\zeta_{y}(E)\big|\Big]^{q}\right)^{1/q}
    \le \mathfrak{c}_q \mathbf{c}(a,A).
\end{equation*} \epr
\subsection{Proof of Lemma \ref{lem:l_p-norm-of-bais}}
First, in view of  the $(\mathfrak{p},\mathfrak{p})$-strong maximal inequality, see e.g. Theorem 9.16 in \cite{WheedenZygmund1977}, one has
\begin{equation*}
%\label{eq1:proof-of-lemma-l_p}
\left\|\overline{\Delta}_{\K,g}(h,\cdot)\right\|_\mathfrak{p}\leq\tau_\mathfrak{p}\left\|\Delta_{\K,g}(h,\cdot)\right\|_\mathfrak{p},
\end{equation*}
where the constant \( \tau_\mathfrak{p} \) depends only of \( \mathfrak{p} \).
Since $\Delta^*_{\K,g}(h,\cdot)\leq \overline{\Delta}_{\K,g}(h,\cdot)+\Delta_{\K,g}(h,\cdot)$  we have
\begin{equation}
\label{eq1:proof-of-lemma-l_p}
\left\|\Delta^*_{\K,g}(h,\cdot)\right\|_\mathfrak{p}\leq(\tau_\mathfrak{p}+1)\left\|\Delta_{\K,g}(h,\cdot)\right\|_\mathfrak{p}.
\end{equation}
 For any $\delta\in (0,h]$ put
$
B(z,\delta)=\Big|\delta^{-1}\int\K\big([u-z]/\delta\big)\big(g(u)-g(z)\big)\dd u\Big|
$
and define
$$
\Delta^{(n)}_{\K,g}(h,z)=\sup_{\delta\in [hn^{-1},h]}B(z,\delta), \;\;n = 1,2, \ldots \;. %\in\mathbb{N}^*.
$$
We remark that the sequence \( \{\Delta^{(n)}_{\K,g}(h,\cdot)\}_{n\geq 1} \) increases monotonically and
\( \Delta^{(n)}_{\K,g}(h,z)\to\Delta_{\K,g}(h,z) \) for any \( z\in\RR \), as \( n\to\infty \).
Hence, by Beppo-Levi's theorem
$$
\left\|\Delta_{\K,g}(h,\cdot)\right\|_\mathfrak{p}=\lim_{n\to\infty}\left\|\Delta^{(n)}_{\K,g}(h,\cdot)\right\|_\mathfrak{p},
$$
and, in view of (\ref{eq1:proof-of-lemma-l_p}), to complete the argument we need to show that
\begin{equation}
\label{eq111:proof-of-lemma-l_p}
\sup_{g\in\mathbb{N}_\mathfrak{p}(s,\cQ)}\left\|\Delta^{(n)}_{\K,g}(h,\cdot)\right\|_\mathfrak{p}\leq 2\cQ h^{s}\|\K\|_\infty
\left[2^{s \mathfrak{p}}-1\right]^{-\frac{1}{\mathfrak{p}}},\;\;\forall n\geq 1.
\end{equation}
Assumption \ref{ass:assumption-on-kernel} ({\it 2}) implies that
we can assert that $B(z,\cdot)$ is continuous on $[n^{-1}h,h]$.
Hence for any $z\in \RR$ there exists $\delta(z)\in [n^{-1}h,h]$ such that
\begin{equation}
\label{eq2:proof-of-lemma-l_p}
\Delta^{(n)}_{\K,g}(h,z)=B\big(z,\delta(z)\big).
%\Delta_{\K,f}(h,z)=\left|\delta^{-1}(z)\int\K\big([u-z]/\delta(z)\big)\big(f(u)-f(z)\big)\dd u\right|.
\end{equation}
For any \( l = 0, \ldots, \log_2n-1 \) (w.l.g. $\log_2n$ is assumed an integer) we consider the slices  $V_l=\big\{z\in\RR:\;\; a_{l+1}<\delta(z)\leq a_l\big\}$ with \( a_l=2^{-l}h \). Later on the integration over empty set is supposed to be zero. Then
\begin{equation}
\label{eq22:proof-of-lemma-l_p}
    \left\| \Delta^{(n)}_{\K,g}(h,\cdot) \right\|^{\mathfrak{p}}_\mathfrak{p}
    =
    \sum_{l=0}^{\log_2n-1} \int_{V_l} |B\big(z,\delta(z)\big)|^{\mathfrak{p}}\dd z.
\end{equation}
We will treat the cases  $s\leq 1$ and $s>1$ separately. If $s<1$ on any slice \( V_l \), \( l = 0, \ldots, \log_2n \),
\begin{eqnarray}
\label{eq33:proof-of-lemma-l_p}
%\hskip-1.0cm
B\big(z,\delta(z)\big)
   &\le&
    \frac{\| \cc K \|_\infty}{\delta(z)}
        \int_{-\frac{\delta(z)}{2} }^{\frac{\delta(z)}{2}}
            \left|g( z + v )-g(z)\right|\dd v
\leq
    \frac{2\| \cc K \|_\infty}{a_l} \int_{-\frac{a_l}{2} }^{\frac{a_l}{2}}
            \left|g( z + v )-g(z)\right|\dd v \nn
&=& 2\|\K\|_\infty \int_{-\frac{1}{2}}^{\frac{1}{2}} \left|g( z + ta_l )-g(z)\right|\dd t.
\end{eqnarray}
We obtain from
(\ref{eq22:proof-of-lemma-l_p}) and (\ref{eq33:proof-of-lemma-l_p}) with the use of Minkowski's  inequality for integrals and writing for ease of notation
$\mu=2\|\K\|_\infty$ that
\begin{eqnarray*}
%&&%\hskip-0.7cm
\left\|\Delta^{(n)}_{\K,g}(h,\cdot)\right\|^{\mathfrak{p}}_\mathfrak{p}
%\sum_{l=0}^\infty\int_{V_l\cap[n^{-1}h,h]}|B\big(z,\delta(z)\big)|^{\mathfrak{p}}\dd z
%\\&&
&\le&
    \mu ^{\mathfrak{p}}
        \sum_{l=0}^{\log_2n-1}
            \int
                 \left|
              \int_{-\frac{1}{2}}^{\frac{1}{2}} \left| g (ta_l+z )-g(z) \right| \dd t
                  \right|^{\mathfrak{p}}\dd z
\nonumber
\\
%&&%\hskip1.6cm
&\le&
    \mu^{\mathfrak{p}}
        \sum_{l=0}^{\log_2n-1}
             \left(
         \int_{-\frac{1}{2}}^{\frac{1}{2}} \left\| g( \cdot + ta_l ) - g(\cdot) \right\|_\mathfrak{p}  \dd t
             \right)^{\mathfrak{p}}
\le \left[\frac{\cQ h^{s }\|\K\|_\infty 2^{1-s}}{(s+1)}\right]^\mathfrak{p}
\sum_{l=0}^\infty 2^{-ls \mathfrak{p}}.
\end{eqnarray*}
Here we have used that $g\in\mathbb{N}_\mathfrak{p}(s,\cQ)$.
Thus, we have for any $s\leq 1$ and any $n\geq 1$
\begin{equation}
\label{eq44:proof-of-lemma-l_p}
\sup_{g\in\mathbb{N}_\mathfrak{p}(s,\cQ)}
\left\|\Delta^{(n)}_{\K,g}(h,\cdot)\right\|_\mathfrak{p}\leq
2\cQ h^{s}\|\K\|_\infty \left[2^{s \mathfrak{p}}-1\right]^{-\frac{1}{\mathfrak{p}}}.
\end{equation}

If $s>1$, using Taylor's formula we have for any $g\in\mathbb{N}_\mathfrak{p}(s,\cQ)$ any $v\in\RR$
\begin{eqnarray*}
g(v+z)-g(z)
=\sum_{m=1}^{m_s}
        \frac{g^{(m)}(z) }{m!}v^m
        +
        \frac{v^{m_s}}{(m_s-1)!}
                \int_0^1 (1-\lambda)^{m_s -1} \left[g^{(m_s)}(z+v\lambda)-g^{(m_s)}(z)\right]\dd\lambda.
\end{eqnarray*}
We have in view of Assumptions \ref{ass:assumption-on-kernel} and \ref{ass2:assumption-on-kernel} for any $z\in\RR$
\begin{equation*}
B\big(z,\delta(z)\big)\leq\frac{\|\K\|_\infty}{(m_s-1)!}\frac{1}{\delta(z)} \int_{-\frac{\delta(z)}{2}}^{\frac{\delta(z)}{2}}
\int_{0}^{1} |v|^{m_s}(1-\lambda)^{m_s -1}\left|g^{(m_s)}\big(z+\lambda v\big)-g^{(m_s)}(z)\right|\dd\lambda\dd v.
\end{equation*}
By the latter inequality for any $z\in V_l$ we get
\begin{equation}
\label{eq3:proof-of-lemma-l_p}
B\big(z,\delta(z)\big)\leq \frac{2\|\K\|_\infty a^{m_s}_l }{(m_s-1)!}
\int_{-\frac{1}{2}}^{\frac{1}{2}}
\int_{0}^{1}|t|^{m_s}(1-\lambda)^{m_s-1}\left|g^{(m_s)}(z + \lambda ta_l )-g^{(m_s)}(z)\right|\dd\lambda\dd t.
\end{equation}
Thus, we obtain from (\ref{eq2:proof-of-lemma-l_p}), (\ref{eq22:proof-of-lemma-l_p}) and (\ref{eq3:proof-of-lemma-l_p}) with the use of Minkowskii inequality for integrals and denoting $\mu=2\|\K\|_\infty \big/(m_s-1)!$ that
\begin{eqnarray*}
&& \left\|\Delta^{(n)}_{\K,f}(h,\cdot)\right\|^{\mathfrak{p}}_\mathfrak{p}
 =
 \sum_{l=0}^{\log_2n-1}  \int_{V_l} |B\big(z,\delta(z)\big)|^{\mathfrak{p}}\dd z
\\
&&\leq
\mu^{\mathfrak{p}}
    \sum_{l=0}^{\log_2n-1}  a^{m_s \mathfrak{p}}_l
        \int
        \left(
        \int_{-\frac{1}{2}}^{\frac{1}{2}}
            \int_{0}^{1}|t|^{m_s}(1-\lambda)^{m_s-1} \left|g^{(m_s)}(z + \lambda ta_l )-g^{(m_s)}(z)\right|\dd\lambda
        \dd t\right)^{\mathfrak{p}}
        \dd z
\nonumber \\
&&\leq
\mu^{\mathfrak{p}}
    \sum_{l=0}^{\log_2n-1} a^{m_s \mathfrak{p}}_l
    \left(
    \int_{-\frac{1}{2}}^{\frac{1}{2}}
        \int_{0}^{1}|t|^{m_s}(1-\lambda)^{m_s-1}\left\|g^{(m_s)}(\cdot + \lambda ta_l )-g^{(m_s)}(\cdot)\right\|_\mathfrak{p}\dd\lambda
        \dd t
    \right)^{\mathfrak{p}}
\nonumber\\
&&\leq \left[\frac{\cQ h^{s }\|\K\|_\infty 2^{1-s}}{(s+1)(m_s+1)(m_s-1)!}\right]^\mathfrak{p}
\sum_{l=0}^\infty 2^{-ls \mathfrak{p}} .
\end{eqnarray*}
Here we have used that $g\in\mathbb{N}_\mathfrak{p}(s,\cQ)$.
Thus, we have for any $s>1$ and $n\geq 1$
\begin{equation}
\label{eq4:proof-of-lemma-l_p}
\sup_{g\in\mathbb{N}_\mathfrak{p}(s,\cQ)
}\left\|\Delta^{(n)}_{\K,g}(h,\cdot)\right\|_\mathfrak{p}\leq 2 \cQ h^{s}
\|\K\|_\infty \left[2^{s \mathfrak{p}}-1\right]^{-\frac{1}{\mathfrak{p}}}.
\end{equation}
We conclude that (\ref{eq111:proof-of-lemma-l_p}) is established in view (\ref{eq44:proof-of-lemma-l_p}) and (\ref{eq4:proof-of-lemma-l_p}).
\epr

%\begin{thebibliography}{9}

\bibliographystyle{agsm}

\begin{thebibliography}{xx}
%
\harvarditem{Bauer et al.}{2009}{BauerHohageMunk}
{\sc  Bauer, F.}, {\sc  Hohage, T.} and {\sc Munk, A.} (2009).
 Iteratively regularized Gauss-Newton method for nonlinear inverse problems with random noise.
 {\it SIAM J. Numer. Anal.} {\bf 47:3} 1827--1846.
 \MR{2505875}.
%
\harvarditem{Bertin and Rivoirard}{2009}{BertinRivoirard}
 {\sc Bertin, K.} and {\sc Rivoirard, V.} (2009).
  Maxiset in sup-norm for kernel estimators.
  {\it TEST} {\bf 18:3} 475--496.
 \MR{2566412}.
%

\harvarditem{Besov  et al.}{1979}{BesovIlNik}
{\sc Besov, O. V.}, {\sc Il'in, V. P.} and {\sc Nikol'skii, S. M.} (1978, 1979). {\it Integral Representations of Functions and Imbedding Theorems.},
 Vol. I,II. Scripta Series in Mathematics., V. H. Winston \& Sons, Washington, D.C.; Halsted Press [John Wiley \& Sons],
  New York-Toronto, Ont.-London.
 \MR{0519341 }, \MR{0521808}.

%
\harvarditem{Chichignoud}{2012}{Chichignoud}
{\sc Chichignoud, M.} (2012).
Minimax and minimax adaptive estimation in multiplicative regression: locally Bayesian approach.
{\em Probab. Theory Related Fields} \textbf{153:3--4} 543--586.
\MR{2948686}.

\harvarditem{Delyon and Juditsky}{1996}{DelyonJuditski1996}
{\sc Delyon, B.} and {\sc Juditsky, A.} (1996). On minimax wavelet estimators.
{\em Appl. Comput. Harmon. Anal.} \textbf{3:3} 215--228.
\MR{1400080}.
%
\harvarditem{Donoho et al.}{1995}{DJKP}
{\sc Donoho, D. L.},  {\sc Johnstone, I. M.},
{\sc Kerkyacharian, G.} and  {\sc Picard, D.}~ (1995).
Wawelet shrinkage: asymptopia? (with discussion).
{\it J. Roy. Statist. Soc.} Ser. R {\bf 57} 301--369.
\MR{1323344}.
%
%
\harvarditem{Ga\"iffas}{2007}{Gaiffas}
{\sc Ga\"iffas, S.} (2007). On pointwise adaptive curve estimation based on inhomogeneous data.
{\em ESAIM P} \&{\em S} 11, 334--364.
\MR{2339297}.
%
%
\harvarditem{Goldenshluger and Lepski}{2008} {GoldLep2008}
{\sc Goldenshluger, A.} and {\sc Lepski, O.} (2008).
Universal pointwise selection rule in multivariate function estimation.
{\em Bernoulli}
{\bf14:4} 1150--1190.
\MR{2543590}.

\harvarditem{Goldenshluger and Lepski}{2009}{GoldLep2009}
{\sc Goldenshluger, A.} and {\sc Lepski, O.}  (2009).
Structural adaptation via \( L_p \)-norm oracle inequalities.
{\it Probab. Theory Related Fields} {\bf 143:1-2} 41--71.
\MR{2449122}.

%
\harvarditem{Ga\"iffas and Lecu\'{e}}{2007}{Lecue}
{\sc Ga\"iffas S.} and {\sc Lecu\'{e} G.} (2007). Optimal rates and adaptation in the single-index model using aggregation.
{\em Electronic J. Statist.}  1, 538--573.
\MR{2369025}.
%


\harvarditem{Golubev}{1992}{Golubev1992}
{\sc Golubev, G. K.} (1992).
%Asymptotic Minimax Estimation of Regression in the Additive Model
Asymptotically minimax estimation of a regression function in an additive model.
{\it Problems Inform. Transmission} {\bf 28:2} 101--112.
\MR{1178413}.

 \harvarditem{Gy\"{o}rfi et al.} {2002} {Gyorfi2002}
 {\sc Gy\"{o}rfi, L.},  {\sc Kohler, M.} {\sc Krzy\.{z}cak, A.} and {\sc Walk, H.} (2002).
 {\it  A distribution-free theory of nonparametric regression.} Springer Series in Statistics. Springer-Verlag, New York. % xvi+647 pp.
\MR{1920390}.

\harvarditem{H\"{a}rdle et al.}{2004} {Haerdle2004}
 {\sc H\"{a}rdle, W.}, {\sc M\"{u}ller, M.}, {\sc Sperlich, S.} and {\sc Werwatz, A.} (2004).
 {\it Nonparametric and semiparametric models.}
 Springer Series in Statistics. Springer-Verlag, New York. %xxviii+299 pp.
 \MR{2061786}.

\harvarditem{Horowitz }{1998}{Horowitz1998}
{\sc Horowitz J. L.} (1998).
{\it Semiparametric and Nonparametric Methods in Econometrics.}
Lecture Notes in Statistics, 131. Springer-Verlag, New York, 1998. %x+204 pp.
\MR{1624936}.

%\bibitem{IbragimovHasm1981}
%Ibragimov, I. A., Has'minskii, R. Z. \emph{Statistical estimation. Asymptotic theory.}
%Translated from the Russian by Samuel Kotz. Applications of Mathematics, 16. Springer-Verlag,
%New York-Berlin, 1981. vii+403 pp.

%
\harvarditem{Kerkyacharian et al.}{2001}{KLP2001}
{\sc Kerkyacharian, G.}, {\sc Lepski, O.} and {\sc Picard, D.}~(2001).
Nonlinear estimation in anisotropic multi--index denoising.
{\em Probab. Theory Related Fields} {\bf 121}, 137--170.
\MR{1863916}.
%
\harvarditem{Kerkyacharian et al.}{2008}{KLP2008}
{\sc Kerkyacharian, G.}, {\sc Lepski, O.} and {\sc Picard, D.}~(2008).
Nonlinear estimation in anisotropic multi--index denoising. Sparce case.
{\em Probab. Theory Appl.} {\bf 52}, 150--171.
\MR{2354574}.
%


\harvarditem{Korostelev and Korosteleva}{2011}{Koros}
{\sc Korostelev, A.} and {\sc Korosteleva, O.} (2011). {\it Mathematical statistics.  Asymptotic minimax theory.} Graduate Studies in Mathematics, 119. American Mathematical Society, Providence, RI.%, 2011. x+246 pp. ISBN: 978-0-8218-5283-5
 \MR{2767163}


%\bibitem{JudLepTsyb2009}
 %Juditsky, A. B., Lepski, O. V. and Tsybakov, A. B. (2009) Nonparametric estimation of composite functions.
 %\emph{Ann. Statist.} \textbf{37:3} 1360--1404.

\harvarditem{Lepski}{1990}{Lep1990}
 {\sc Lepskii, O. V.} (1990).
 A problem of adaptive estimation in Gaussian white noise.
 %(Russian) \textit{Teor. Veroyatnost. i Primenen.} \textbf{35:3} 459--470;  translation in
{\it Theory Probab. Appl.} {\bf 35:3} 454--466.
\MR{1091202}.


 \harvarditem{Lepski and Levit}{1999}{LepskiLevit99}
{\sc Lepski, O. V.} and {\sc Levit, B. Y.} (1999).  Adaptive nonparametric estimation of smooth multivariate functions.
{\it Math. Methods Statist.} {\bf 8:3} 344--370.
\MR{1735470}.

\harvarditem{Lepski et al.}{1997}{LepMamSpok97}
{\sc Lepski, O. V.}, {\sc Mammen, E.} and \textsc{Spokoiny, V.G.} (1997).
Optimal spatial adaptation to inhomogeneous smoothness: an approach based on kernel estimates with variable bandwidth selectors.
{\it Ann. Statist.} {\bf 25:3} 929--947.
\MR{1447734}.

%\bibitem{LepSpok97}
%Lepski, O. V.,  and Spokoiny, V.G. (1997).
%Optimal pointwise adaptive methods in nonparametric estimation.
%\textit{Ann. Statist.} \textbf{25:6} 2512--2546.


\harvarditem{Lifshits}{1995}{Lif}
\textsc{Lifshits M.A.} (1995).
\textit{Gaussian Random Functions}.
Mathematics and its Applications, 322. Kluwer Academic Publishers, Dordrecht. %xii+333 pp.
\MR{1472736}.

\begin{comment}
\harvarditem{Maddala}{1983}{MADDALA}
\textsc{Maddala, G.} (1983).
\textit{Limited-Dependent and Qualitative Variables in Econometrics.}
 Econometric Society Monographs in Quantitative Economics, 3. Cambridge University Press,
 Cambridge, 1983. %xii+401 pp.
\MR{0799154}
\end{comment}
%
%
\harvarditem{Serdyukova}{2012}{Serdyukova}
{\sc Serdyukova, N.} (2012).
Spatial adaptation in heteroscedastic regression: propagation approach.
{\it Electron. J. Stat.} {\bf 6} 861--907.
\MR{2988432 }
%
\harvarditem{Stone} {1985} {Stone85}
\textsc{Stone, C.J.} (1985).
Additive regression and other nonparametric models. {\it Ann. Statist.} \textbf{13:2} 689--705.
\MR{0790566}

%\bibitem{Talagrand}
%{ Talagrand, M.}~(1994).
%Sharper bounds for Gaussian and empirical processes.
%{\em Ann. Probab.} {\bf 22}, 28--76.

\harvarditem{Wheeden and Zygmund}{1977}{WheedenZygmund1977}
\textsc{Wheeden, R. L.} and \textsc{Zygmund, A.} (1977).
\textit{Measure and integral. An introduction to real analysis. }
Pure and Applied Mathematics, Vol. 43.
Marcel Dekker, Inc., New York-Basel. %x+274 pp.
\MR{0492146}.

\end{thebibliography}

\end{document}